\documentclass[a4paper]{elsarticle}

\usepackage[utf8]{inputenc}
\usepackage{import}

\usepackage{fontawesome}
\usepackage{lipsum}

\usepackage[english]{babel}

\usepackage[margin=2.45cm]{geometry}

\usepackage{xcolor}

\usepackage{graphicx}
\usepackage{svg}
\usepackage{epstopdf}
\usepackage{tikz}
\usetikzlibrary{arrows.meta,calc,positioning}
\usepackage{pgfplots}
\pgfplotsset{compat=newest}
\usepackage[percent]{overpic}
\usepackage{placeins}

\usepackage{circledsteps}
\pgfkeys{/csteps/fill color=white}
\pgfkeys{/csteps/inner ysep=5pt}
\pgfkeys{/csteps/inner xsep=5pt}

\usepackage{enumitem}
\setlist[1]{itemsep=-5pt}

\usepackage{amsmath}
\usepackage{amssymb,amsthm}
\usepackage{cases}
\usepackage{mathrsfs}
\usepackage{euscript}
\usepackage{bm}
\usepackage{stmaryrd}
\usepackage{mathtools}
\usepackage{empheq}

\theoremstyle{definition}

\usepackage{algorithm}
\usepackage{algpseudocode}

\usepackage{appendix}

\usepackage{pdflscape}
\usepackage{setspace}
\usepackage{framed}
\usepackage{pdfpages}

\usepackage{tabularx}
\usepackage{booktabs}

\usepackage{multicol}
\usepackage{multirow}
\usepackage{xparse}
\usepackage{xspace}
\usepackage{makecell}
\usepackage{tablefootnote}

\usepackage{footnote}

\usepackage{caption}
\usepackage{subcaption}

\usepackage[notref,notcite,final]{showkeys}

\usepackage[noprefix]{nomencl}
\makenomenclature

\renewcommand*\nompreamble{\begin{multicols}{2}}
\renewcommand*\nompostamble{\end{multicols}}

\usepackage{etoolbox}
\usepackage{stackengine}
\setstackgap{S}{1.75pt}
\usepackage{comment}

\biboptions{sort&compress}
\usepackage[numbers]{natbib}
\usepackage[colorlinks=true]{hyperref}
\usepackage{cleveref}

\def \usetodonotes{} 

\ifdefined\usetodonotes

\definecolor{C0}{HTML}{1F77B4}
\definecolor{C1}{HTML}{FF7F0E}
\definecolor{C2}{HTML}{2ca02c}
\definecolor{C3}{HTML}{d62728}
\definecolor{C4}{HTML}{9467bd}
\definecolor{C5}{HTML}{8c564b}
\definecolor{TileBlue}{HTML}{78BFDC}
\definecolor{TileOrange}{HTML}{F8AC6B}

\usepackage[textsize=scriptsize]{todonotes}

\NewDocumentCommand{\JU}{O{} m}{%
    \todo[fancyline,color=C0!10,linecolor=C0,bordercolor=C0,#1]{JU: #2}%
}
\NewDocumentCommand{\PR}{O{} m}{%
    \todo[fancyline,color=C1!10,linecolor=C1,bordercolor=C1,#1]{PR: #2}%
}

\newlength{\todostripwidth}
\usepackage{atbegshi}
\AtBeginShipout{%
    \AtBeginShipoutUpperLeft{%
        \begin{tikzpicture}[overlay,remember picture]
            \fill [gray!0] (current page.north east) rectangle
                ($ (current page.south east) - (\todostripwidth,0cm) $);
        \end{tikzpicture}
    }
}

\else

\NewDocumentCommand{\JU}{O{} m}{}
\NewDocumentCommand{\PR}{O{} m}{}

\fi

\newcommand{\md}{\kern0.1em---\kern0.1em\ignorespaces}
\newcommand{\nd}{\kern0.1em--\kern0.1em\ignorespaces}

\newcommand{\RR}{\bm Q}
\newcommand{\Cper}{\bm{\mathsf C}_{\mathrm{per}}}

\makeatletter
\newcommand{\multiline}[1]{%
  \begin{tabularx}{\dimexpr\linewidth-\ALG@thistlm}[t]{@{}X@{}}
    #1
  \end{tabularx}
}
\makeatother

\definecolor{myblue}{RGB}{0, 0, 204}

\allowdisplaybreaks
\graphicspath{
  {./}
  {./source/}
}
\begin{document}

\begin{frontmatter}

\author{Weiyuan Fan}
\ead{weiyuanf@umich.edu}
\author{Jacinto Ulloa\corref{cor1}}
\ead{julloa@umich.edu}

\address{Department of Mechanical Engineering, University of Michigan, Ann Arbor, MI 48109, USA}

\cortext[cor1]{Corresponding author}

\title{Higher-dimensional homogenization of quasiperiodic elastic materials}

\begin{abstract}
Quasiperiodic geometrical arrangements, which exhibit long-range order in the absence of periodicity, have attracted interest in the design of composites and mechanical metamaterials. In this context, direct numerical simulations explicitly resolving the quasiperiodic microstructure become prohibitively expensive for macroscopic problems, motivating the development of homogenized continuum models. However, the absence of a periodic unit cell in physical space hinders the use of conventional homogenization methods. In the present work, we exploit higher-dimensional periodic representations to perform computational homogenization of quasiperiodic elastic composites. Thereby, effective properties are obtained from unit-cell problems solved using the fast Fourier transform in higher-dimensional spaces. Using two-dimensional Penrose tilings and three-dimensional Ammann--Kramer--Neri tilings as paradigmatic examples, we describe the distribution of heterogeneous material properties in the corresponding higher-dimensional unit cells, compute effective elasticity tensors, and assess convergence of fully resolved finite patches to the homogenized limits. We provide open-source implementations of all numerical models.   
\end{abstract}


\end{frontmatter}
\setcounter{tocdepth}{2}
\tableofcontents

\section{Introduction}
Since their discovery in metallic alloys and subsequent realization in architected media, quasicrystals have attracted significant interest across solid mechanics and materials physics~\citep{shechtman1984metallic, levine1984quasicrystals, montemayor2015materials, penrose1974role, kramer1984periodic, lu2022architectural}. In particular, quasiperiodic composites and metamaterials\md such as two-dimensional Penrose tilings with 5-fold rotational symmetry and three-dimensional Ammann--Kramer--Neri (AKN) tilings with icosahedral symmetry\md exhibit deterministic long-range order without periodicity (or translational symmetry), providing an alternative to both periodic and random structures. These structures have been shown to yield,  beyond their aesthetic appeal, macroscopically isotropic responses with enhanced stiffness, deformability, and wave-attenuation properties~\citep{wang2020quasiperiodic, rosa2025enhanced, beli2021mechanics, beli2022wave}. Harnessing these properties in macrostructural applications requires predictive models that bridge the microscopic and macroscopic scales, which poses a pressing computational challenge due to the lack of periodicity of the underlying microstructure.

Direct numerical simulations are computationally intractable for composites or metamaterials with a marked scale separation between individual tile dimensions and structural boundaries. Standard multiscale approaches~\cite{kochmann2019multiscale} address this problem by evaluating effective constitutive properties over a unit cell subject to periodic boundary conditions, as commonly employed for periodic composites~\cite{michel1999effective, cioranescu1999introduction, nemat2011homogenization, Milton2002, andreassen2014determine} and mechanical lattices~\cite{Fleck2010,Vigliotti2012,Phlipot2019,Ariza2024,Ulloa2025}. Quasicrystals, however, admit no finite translationally periodic unit cell in the physical space ($\mathbb{R}^2$ or $\mathbb{R}^3$). Truncating physical domains to enforce artificial periodicity introduces boundary-layer distortions, artificial size dependence, and sample-to-sample variance from finite patches. On the other hand, stochastic homogenization can identify deterministic effective properties for disordered media under suitable stationarity and ergodicity assumptions~\cite{Milton2002,Torquato2002,gloria2015quantification}, but does not explicitly exploit the deterministic algebraic structure governing quasicrystalline geometry.

While quasiperiodic geometries lack translational symmetry in physical space, they admit periodic representations in higher-dimensional spaces, constituting a form of \emph{deterministic disorder}. This higher-dimensional periodicity has been exploited in two-scale cut-and-projection methods to formulate higher-dimensional unit-cell problems for quasiperiodic homogenization, without constructing periodic approximations of the physical microstructure. In this setting, quasicrystals in the physical space $\mathbb R^n$ are represented by periodic structures in $\mathbb R^m$, with $m>n$. While rigorous mathematical foundations have been established in this context~\citep{braides2009,bouchitte2010homogenization,Wellander2019,ferreira2021homogenization}, to our knowledge, the underlying theory has yet to be exploited for the computational homogenization of quasiperiodic elastic composites.

In the present work, we address this gap by developing a first computational framework for the homogenization of quasiperiodic elastic composites and assessing their convergence to the resulting homogenized limits. The effective models are derived from a variational standpoint, with the effective elastic properties computed from periodic unit-cell problems solved using the fast Fourier transform (FFT) in higher-dimensional spaces. We consider, as paradigmatic examples,  Penrose tilings ($n=2$, $m=4$) and AKN icosahedral tilings ($n=3$, $m=6$). For both cases, we describe the geometrical construction of heterogeneous material properties in the corresponding higher-dimensional periodic unit cell. We then compute effective properties and assess convergence by comparing direct numerical solutions of finite quasiperiodic patches under loading with their homogenized counterparts. We provide open-source computational implementations for both examples, including the higher-dimensional spectral solver in \texttt{MATLAB} and the physical-space finite element models implemented in \texttt{Python} using \texttt{FEniCS}.


This paper is structured as follows. Section~\ref{sec:geometry} describes the geometry and material distributions of the quasiperiodic solids, setting the stage for the homogenization framework presented in section~\ref{sec:hom}. Section~\ref{sec:comp_imp} details the construction of the higher-dimensional periodic cells and the FFT solution of the unit-cell problem. Section~\ref{sec:examples} presents the numerical simulations and convergence studies, while section~\ref{sec:conclusions} closes the paper with conclusions and perspectives. 

\section{Geometrical description of quasiperiodic materials}
\label{sec:geometry}

We consider two quasiperiodic architectures: the two-dimensional Penrose tiling and the three-dimensional Ammann--Kramer--Neri (AKN) tiling. Both are generated through multigrid constructions and admit periodic representations in higher-dimensional spaces. These representations allow the heterogeneous elasticity tensors to be expressed as periodic tensors, providing the geometrical foundation for quasiperiodic homogenization.

\subsection{Two-dimensional Penrose tiling}

\subsubsection{Two-dimensional pentagrid construction}\label{sec:element_label}

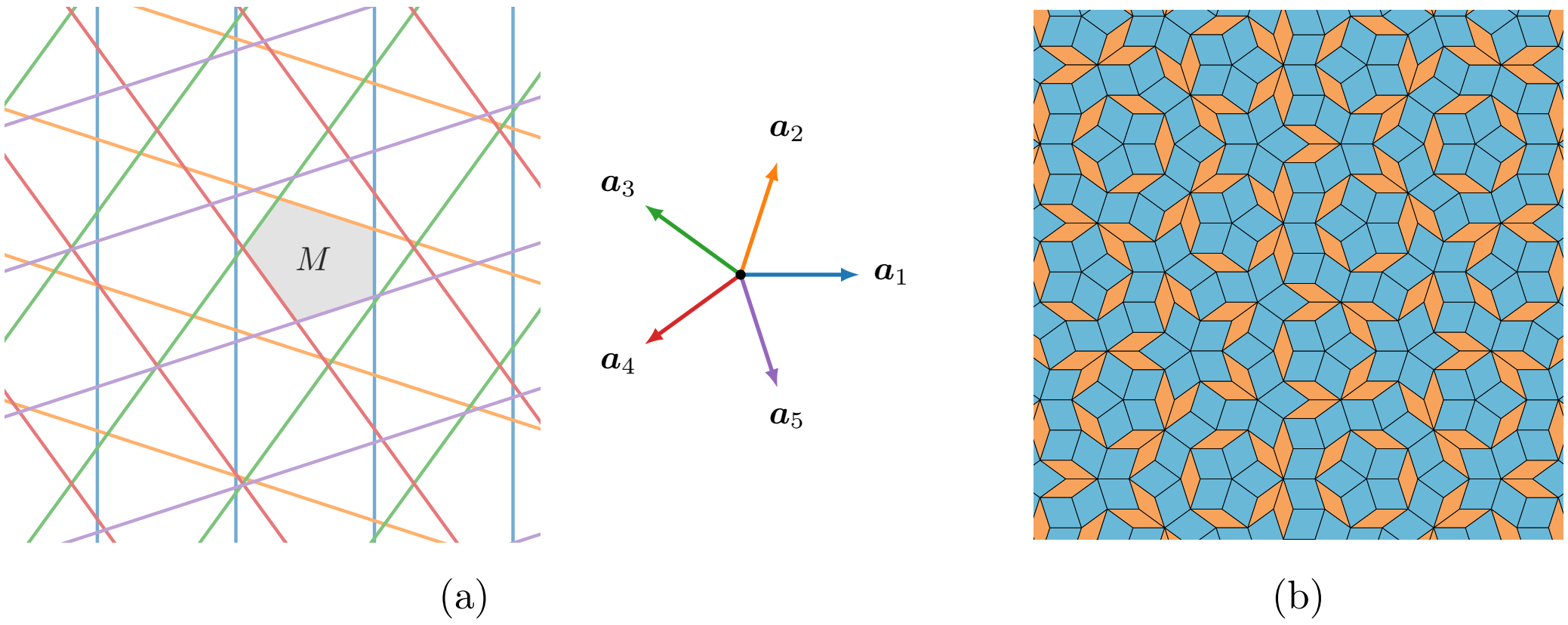
\begin{figure}[htbp]
    \centering
    \begin{minipage}[c]{0.5\textwidth}
        \centering
            \resizebox{\linewidth}{!}{%
                \begin{tikzpicture}[>=Latex]
                    \begin{scope}[xshift=-2.0cm]
                        \clip (-2.80,-2.80) rectangle (2.80,2.80);
                        \begin{scope}[shift={(0.43,0.17)},scale=1.45,shift={(-0.43,-0.17)}]
                            \fill[gray!22]
                                (0.87000,-0.10636) --
                                (0.87000,0.37974) --
                                (0.17936,0.60414) --
                                (-0.08034,0.24669) --
                                (0.30867,-0.28875) -- cycle;

                            \foreach \j/\g/\col in {1/0.13/C0,2/0.37/C1,3/-0.21/C2,4/0.08/C3,5/-0.37/C4}{
                                \pgfmathsetmacro{\ang}{72*(\j-1)}
                                \foreach \k in {-7,-6,...,7}{
                                    \draw[draw=\col!62,line width=1pt]
                                        ({(\k-\g)*cos(\ang)-8*sin(\ang)},
                                         {(\k-\g)*sin(\ang)+8*cos(\ang)})
                                        --
                                        ({(\k-\g)*cos(\ang)+8*sin(\ang)},
                                         {(\k-\g)*sin(\ang)-8*cos(\ang)});
                                }
                            }
                            \node[text=black!85] at (0.43,0.17) {$M$};
                        \end{scope}
                    \end{scope}

                    \begin{scope}[xshift=2.9cm]
                        \foreach \j/\col in {1/C0,2/C1,3/C2,4/C3,5/C4}{
                            \pgfmathsetmacro{\ang}{72*(\j-1)}
                            \draw[draw=\col,very thick,-{Latex[length=2mm]}]
                                (0,0)--({1.25*cos(\ang)},{1.25*sin(\ang)});
                            \node[fill=white,inner sep=0.8pt]
                                at ({1.58*cos(\ang)},{1.58*sin(\ang)})
                                {$\bm a_{\j}$};
                        }
                        \fill (0,0) circle (1.5pt);
                    \end{scope}
                \end{tikzpicture}
            }
    \end{minipage}%
    \hspace{0.06\textwidth}%
    \begin{minipage}[c]{0.288\textwidth}
        \centering
        \begin{overpic}[width=\linewidth]
            {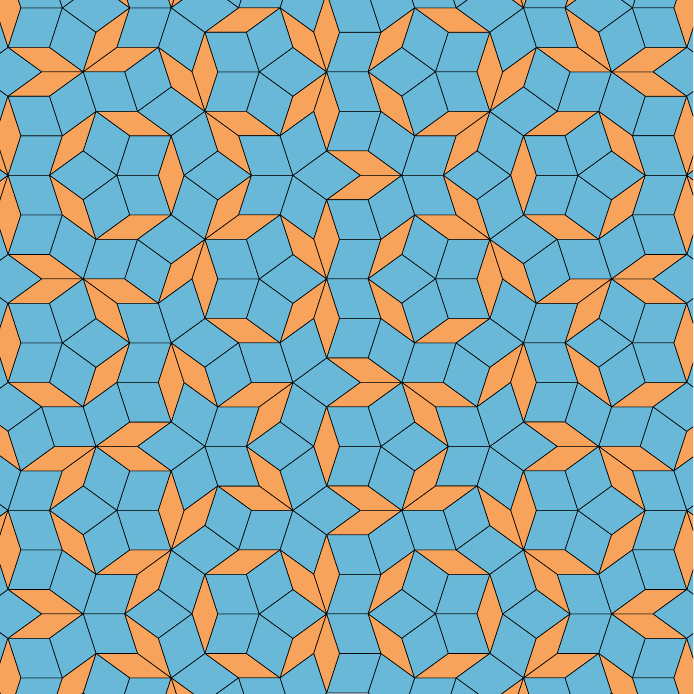}
        \end{overpic}
    \end{minipage}
    \par\vspace{0.8em}
        \noindent
        \begin{minipage}{0.5\textwidth}
            \centering
            (a)
        \end{minipage}%
        \hspace{0.06\textwidth}%
        \begin{minipage}{0.288\textwidth}
            \centering
            (b)
        \end{minipage}
    \caption{Pentagrid and corresponding Penrose tiling. (a) Regular pentagrid for $\bm\gamma = (0.13,0.37,-0.21,0.08,-0.37)^\mathsf{T}$, with an element $M$ highlighted in gray and the unit vectors $\{\bm a_j\}_{j=1}^5$ normal to the correspondingly colored families of parallel lines. (b) Finite Penrose patch, with thick and thin rhombi shown in blue and orange, respectively.}
    \label{fig:pentagrid_element}
\end{figure}

We describe the geometric construction of a two-dimensional Penrose tiling following the pentagrid method of \citet{DeBruijn1981}. Consider five families of parallel lines whose unit-normal directions are equally spaced at angular intervals of $72^\circ$. The corresponding unit vectors are
\begin{equation}\label{eq:five_directions}
    \bm a_j=\begin{pmatrix}\cos(2\pi(j-1)/5)\\\sin(2\pi(j-1)/5)\end{pmatrix},\qquad j=1,\ldots,5.
\end{equation}
Let us introduce the linear map $\bm A^5:\mathbb R^2\to\mathbb R^5$, represented by $\bm A^5=(\bm a_1, \dots, \bm a_5)^\mathsf T\in\mathbb R^{5\times2}$, and an offset vector $\bm\gamma\in\mathbb R^5$, satisfying $\sum_{j=1}^5\gamma_j=0$. Furthermore, for $\bm z\in\mathbb R^2$, consider a vector $\bm q(\bm z)\in\mathbb R^5$ of shifted projections,
\begin{equation}\label{eq:five_dimensional_point}
    \bm q(\bm z)=\bm\gamma+\bm A^5\bm z.
\end{equation}
For every $k\in\mathbb Z$, the $k$th line of family $j$ is then given by
\begin{equation}
    \label{eq:grid_lines}
    G_{j}^k
    =
    \big\{
        \bm z\in\mathbb R^2 \ : \
        q_j(\bm z)=k
    \big\}.
\end{equation}
For fixed $j$, the collection $\{G_{j}^k\}_{k\in\mathbb Z}$ represents a single family of parallel lines separated by unit distance. The union of the five families is called the \emph{pentagrid}.
The pentagrid lines partition $\mathbb R^2$ into polygonal regions, each of which is referred to as a \emph{pentagrid element} $M$. Figure~\ref{fig:pentagrid_element}(a) highlights an element within a pentagrid and shows the five unit normals $\{\bm a_j\}_{j=1}^5$ defining the families of parallel lines. 

\begin{figure}[b!]
    \centering
    \begin{tikzpicture}[>=Latex]
        \begin{scope}[xshift=-4.25cm]
            \fill[gray!12] (-2.0,-1.45) -- (0,-1.45) -- (0,0) -- (-2.0,{2*tan(18)}) -- cycle;
            \draw[C0,thick] (0,-1.45) -- (0,1.45);
            \draw[C2,thick] (-2.0,{2*tan(18)}) -- (2.0,{-2*tan(18)});
            \node at (-1.18,-0.36) {$\bm\kappa$};
            \node at (0.75,-0.72) {$\bm\kappa+\bm g_r$};
            \node at (-0.75,0.75) {$\bm\kappa+\bm g_s$};
            \node at (1.1,0.5) {$\bm\kappa+\bm g_r+\bm g_s$};
            \draw[C0,very thick,-{Latex[length=2mm]}]
                (-0.58,-1.18)--(0.58,-1.18);
            \node[C0] at (-0.4,-0.9) {$+\bm g_r$};
            \fill (0,0) circle (1.6pt);
            \draw[C0,thick] (-1.85,-1.95)--(-1.52,-1.95);
            \node[anchor=west] at (-1.42,-1.95) {family $r$};
            \draw[C2,thick] (0.40,-1.95)--(0.73,-1.95);
            \node[anchor=west] at (0.83,-1.95) {family $s$};
            \node at (0,-2.45) {pentagrid};
        \end{scope}
        \draw[->,very thick] (-1,-0.33) -- (1.30,-0.33)
            node[midway,above] {\!\!\!\!dual tiling};
        \begin{scope}[xshift=2.4cm,yshift=-1.33cm]
            \coordinate (P) at (0,0);
            \coordinate (R) at (2.15,0);
            \coordinate (S) at ({2.15*cos(72)},{2.15*sin(72)});
            \coordinate (RS) at ({2.15+2.15*cos(72)},{2.15*sin(72)});
            \fill[C1!18] (P)--(R)--(RS)--(S)--cycle;
            \draw[C1,very thick] (P)--(R)--(RS)--(S)--cycle;
            \draw[C0,very thick,-{Latex[length=2mm]}]
                (P)--(R) node[midway,below=2pt] {$\bm a_r$};
            \fill (P) circle (1.5pt) node[below left] {$\bm p$};
            \fill (R) circle (1.5pt) node[below right] {$\bm p+\bm a_r$};
            \fill (S) circle (1.5pt) node[above left] {$\bm p+\bm a_s$};
            \fill (RS) circle (1.5pt)
                node[above right] {$\bm p+\bm a_r+\bm a_s$};
            \node[below=0.8cm] at (1.38,0) {Penrose rhombus};
        \end{scope}
    \end{tikzpicture}
    \caption{Duality between crossings of pentagrid line families $r$ and $s$ and a unit-edge Penrose rhombus. Crossing a line of family $j\in\{r,s\}$ changes the element label $\bm\kappa\in\mathbb{Z}^5$ by $\pm\bm g_j$ and the corresponding physical vertex $\bm p\in\mathbb{R}^2$ by $\pm\bm a_j$.}
    \label{fig:crossing_rhombus}
\end{figure}
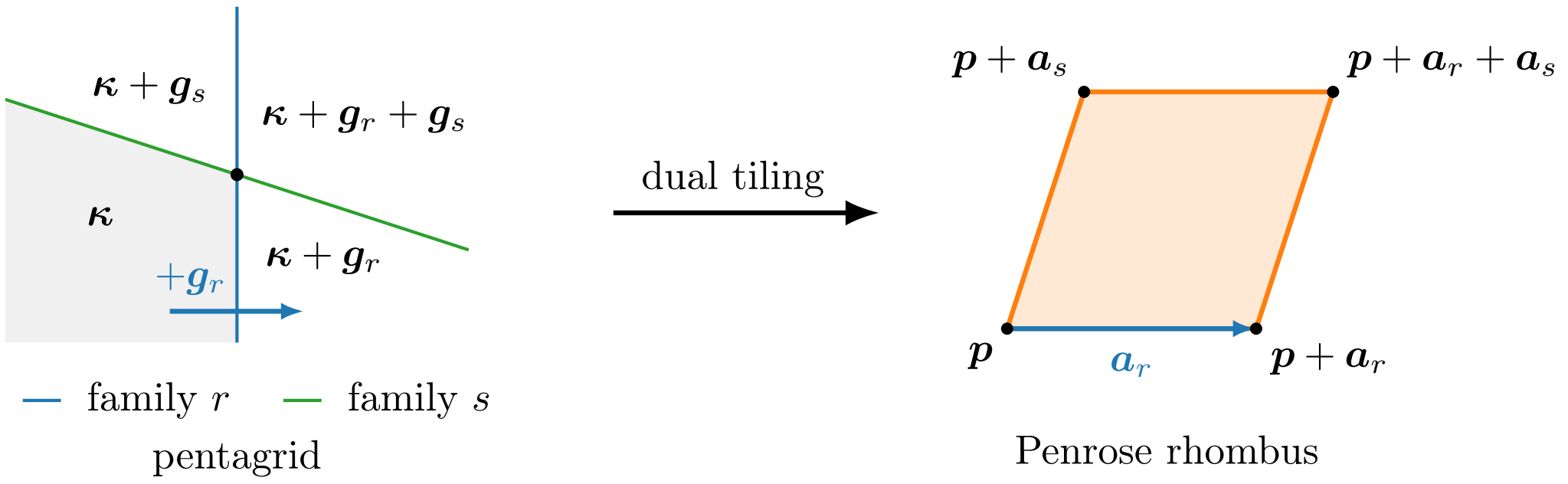

\subsubsection{Dual tiling}
\label{sec:crossing_to_rhombus}

Assuming a regular pentagrid, meaning that no point lies on three or more lines, each pentagrid element $M$ is assigned an integer label
\begin{equation}\label{eq:penrose_label}
\bm\kappa(M)=\big(\lceil q_1(\bm z)\rceil,\ldots,\lceil q_5(\bm z)\rceil\big)^\mathsf T\in\mathbb Z^5,\qquad\bm z\in M.
\end{equation}
The label is independent of the choice of $\bm z\in M$. 

The \emph{dual} tiling is constructed by the mapping $M\mapsto\bm p(M)\in\mathbb R^2$ from pentagrid elements to vertices in the physical plane. Specifically,
\begin{equation}\label{eq:dual_vertex}
    \bm p(M)=\sum_{j=1}^5\kappa_j(M)\bm a_j=(\bm A^5)^\mathsf T\bm\kappa(M).
\end{equation}
Let $M$ and $M'$ share a boundary segment belonging to family $j$.  Their labels satisfy
\begin{equation}
    \label{eq:adjacent_element_labels}
    \bm\kappa(M')-\bm\kappa(M)=\pm\bm g_j,
\end{equation}
where $\bm g_j$ denotes the $j$th standard basis vector of $\mathbb R^5$. Thus, crossing a line of family $j$ changes only component $j$ of the pentagrid element label by $\pm 1$, where the sign depends on the direction in which the line is crossed. Subtracting equation~\eqref{eq:dual_vertex} for the two pentagrid elements gives
\begin{equation}
    \label{eq:dual_edge}
    \bm p(M')-\bm p(M)
    =
    \sum_{i=1}^5
    \bigl[\kappa_i(M')-\kappa_i(M)\bigr]\bm a_i
    =
    \pm\bm a_j.
\end{equation}
Note that all terms except the term with $i=j$ vanish. The segment with endpoints $\bm p(M)$ and $\bm p(M')$ is the Penrose edge dual to the common pentagrid segment. Therefore, pentagrid elements correspond to Penrose vertices and pentagrid segments correspond to Penrose edges.

Suppose that a line of family $r$ crosses a line of family $s$, with $r\neq s$. Regularity implies that four pentagrid elements meet at the crossing, generating a single Penrose tile (Figure~\ref{fig:crossing_rhombus}). Specifically, denoting one of the labels by $\bm\kappa$ and its corresponding physical point by $\bm p$, the four labels are given by
\begin{equation*}
    \bm\kappa,
    \qquad
    \bm\kappa+\bm g_r,
    \qquad
    \bm\kappa+\bm g_s,
    \qquad
    \bm\kappa+\bm g_r+\bm g_s.
\end{equation*}
Using equation~\eqref{eq:dual_vertex}, these labels give the physical points
\begin{equation*}
    \bm p,
    \qquad
    \bm p+\bm a_r,
    \qquad
    \bm p+\bm a_s,
    \qquad
    \bm p+\bm a_r+\bm a_s,
\end{equation*}
which correspond to the four vertices of a rhombus, with $\lvert\bm a_r\rvert=\lvert\bm a_s\rvert=1$. The acute interior angle of the resulting rhombus is either $36^\circ$ or $72^\circ$, producing, respectively, the \emph{thin} and \emph{thick} Penrose tiles.

\begin{figure}[t!]
    \centering
    \begin{tikzpicture}[>=Latex]
        \begin{scope}[xshift=-3.8cm]
            \fill[C1!12] (0,1)--(-1.90,{1+1.90*tan(18)})--(-1.90,{-1-1.90*tan(18)})--(0,-1)--cycle;
            \fill[gray!10] (0,1)--(1.90,{1-1.90*tan(18)})--(1.90,{-1+1.90*tan(18)})--(0,-1)--cycle;
            \draw[C0,thick] (0,-1.55)--(0,1.55);
            \draw[C1,thick] (-1.90,{1+1.90*tan(18)})--(1.90,{1-1.90*tan(18)});
            \draw[C2,thick] (-1.90,{-1-1.90*tan(18)})--(1.90,{-1+1.90*tan(18)});
            \draw[black,very thick] (0,-1)--(0,1);
            \fill[C1] (0,1) circle (1.8pt);
            \node[above left=4pt] at (0.05,1.1) {$C_1$};
            \fill[C2] (0,-1) circle (1.8pt);
            \node[below left=4pt] at (0.05,-1.1) {$C_2$};
            \node at (-0.72,0) {$M$};
            \node at (0.72,0) {$M'$};
            \node at (0,-2.25) {pentagrid};
        \end{scope}
        \draw[->,very thick] (-0.9,0)--(1.4,0)
            node[midway,above] {\!\!\!\!dual tiling};
        \begin{scope}[xshift=3.75cm]
            \coordinate (P) at (-0.75,0);
            \coordinate (Q) at (0.75,0);
            \coordinate (U) at ({-0.75+1.5*cos(72)},{1.5*sin(72)});
            \coordinate (V) at ({0.75+1.5*cos(72)},{1.5*sin(72)});
            \coordinate (L) at ({-0.75+1.5*cos(72)},{-1.5*sin(72)});
            \coordinate (W) at ({0.75+1.5*cos(72)},{-1.5*sin(72)});
            \fill[C1!18] (P)--(Q)--(V)--(U)--cycle;
            \fill[C2!18] (P)--(Q)--(W)--(L)--cycle;
            \draw[C1,thick] (P)--(Q)--(V)--(U)--cycle;
            \draw[C2,thick] (P)--(Q)--(W)--(L)--cycle;
            \draw[black,very thick] (P)--(Q);
            \fill (P) circle (1.5pt)
                node[left=2pt] {$\bm p(M)$};
            \fill (Q) circle (1.5pt)
                node[right=2pt] {$\bm p(M')$};
            \node at (0.23,0.72) {$C_1$};
            \node at (0.23,-0.72) {$C_2$};
            \node at (0.23,-2.25) {Penrose tiling};
        \end{scope}
    \end{tikzpicture}
    \caption{Duality between adjacent pentagrid elements and Penrose vertices, where the consecutive crossings $C_1$ and $C_2$ joined by a pentagrid segment (left) correspond to Penrose rhombi sharing the edge $[\bm p(M),\bm p(M')]$ (right).}
    \label{fig:adjacent_rhombi}
\end{figure}
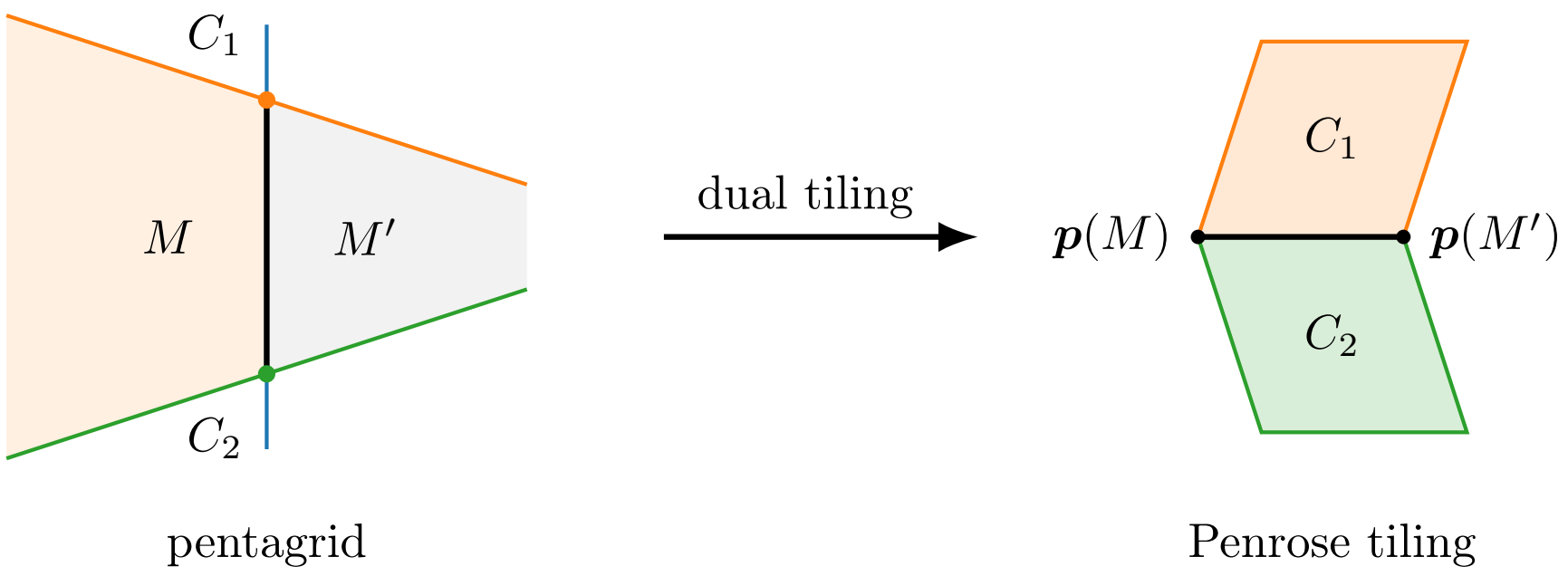

To see how these rhombi assemble, consider two consecutive crossings on a pentagrid line segment, with $M$ and $M'$ the neighboring elements (Figure~\ref{fig:adjacent_rhombi}). The two crossings produce two rhombi that share the edge $[\bm p(M),\bm p(M')]$ and are thus adjacent in the physical plane. Accordingly, applying this construction to all pentagrid crossings produces rhombi that cover the physical plane without gaps or overlapping interiors~\cite{DeBruijn1981}. Figure~\ref{fig:pentagrid_element}(b) shows a finite Penrose patch constructed in this way.

\subsubsection{Four-dimensional periodic representation}
\label{sec:five_dimensional_lattice}

As $\bm z$ varies over $\mathbb R^2$, the points $\bm q(\bm z)\in\mathbb R^5$ given by equation~\eqref{eq:five_dimensional_point} form the shifted plane~\citep{braides2009}
\begin{equation}
	\label{eq:pi_plane}
    \Pi_{\bm\gamma}^5
    =
    \big\{
        \bm q(\bm z) \ : \
        \bm z\in\mathbb R^2
    \big\}
    =
    \bm\gamma+\bm A^5\mathbb R^2.
\end{equation}
Moreover, we associate with each $\bm\kappa\in\mathbb Z^5$ the open unit hypercube
\begin{equation}
    \label{eq:label_cube}
    Q^5_{\bm\kappa}
    =
    \big\{
        \bm q\in\mathbb R^5 \ : \
        \kappa_j-1<q_j<\kappa_j,
        \ j=1,\ldots,5
    \big\}.
\end{equation}
The hypercube closures form the standard $\mathbb Z^5$-periodic grid in $\mathbb R^5$, with neighboring hypercubes sharing a four-dimensional face. The shifted plane $\Pi_{\bm\gamma}^5$ passes through this periodic grid, and equations~\eqref{eq:penrose_label}, \eqref{eq:pi_plane}, and \eqref{eq:label_cube} give, for each pentagrid element $M$,
\begin{equation}
    \big\{
        \bm q(\bm z) \ : \
        \bm z\in M
    \big\}
    =
    \Pi_{\bm\gamma}^5\cap Q^5_{\bm\kappa(M)}.
\end{equation}
Thus, under the mapping~\eqref{eq:five_dimensional_point}, each $M$ is identified uniquely with a two-dimensional section $\Pi_{\bm\gamma}^5\cap Q^5_{\bm\kappa(M)}$.

Since $\sum_{j=1}^5\bm a_j=\bm0$ and $\sum_{j=1}^5\gamma_j=0$, the shifted plane~\eqref{eq:pi_plane} lies in a four-dimensional subspace $\Xi$, i.e.,
\begin{equation}
    \label{eq:Xi^4_definition}
    \Pi_{\bm\gamma}^5\subset \Xi, \qquad \Xi
    =
    \bigg\{
        \bm q\in\mathbb R^5 \ : \
        \sum_{j=1}^5 q_j = 0
    \bigg\}.
\end{equation}
Note that $\Xi\cong\mathbb R^4$ is a four-dimensional hyperplane in $\mathbb R^5$. For each $\bm\kappa\in\mathbb Z^5$, let $S_{\bm\kappa}^5=\Xi\cap Q^5_{\bm\kappa}$, such that $\Pi_{\bm\gamma}^5\cap Q^5_{\bm\kappa(M)}
    \subset
    S_{\bm\kappa(M)}^5$. The nonempty sets $\{S_{\bm\kappa}^5\}_{\bm\kappa\in\mathbb Z^5}$ form a partition of $\Xi$, inheriting the periodicity of the five-dimensional hypercube grid under integer translations satisfying $\sum_{j=1}^5n_j=0$, i.e., $S_{\bm\kappa+\bm n}^5 = S_{\bm\kappa}^5+\bm n$. These translations form a four-dimensional lattice, with a possible choice of basis~(cf.~\citet{Mazac2023})
\begin{equation}
    \label{eq:b_basis}
    \bm b_\alpha=\bm g_\alpha-\bm g_5,\qquad \alpha=1,\ldots,4.
\end{equation}
Assembling these vectors into $\bm B=(\bm b_1\ \bm b_2\ \bm b_3\ \bm b_4)\in\mathbb R^{5\times4}$, we have $\Xi=\bm B\mathbb R^4$, and every point $\bm w\in\Xi$ can be uniquely represented as $\bm w=\bm B\bm y$, with $\bm y\in\mathbb R^4$. Moreover, we have $S_{\bm\kappa}^5=\bm B S_{\bm\kappa}^4$, with
\begin{equation}
    \label{eq:four_dimensional_partition}
    S_{\bm\kappa}^4=\left\{\bm y\in\mathbb R^4 \ : \ \bm B\bm y\in S_{\bm\kappa}^5\right\}.
\end{equation}
Since $S_{\bm\kappa+\bm B\bm m}^4=S_{\bm\kappa}^4+\bm m$ for $\bm m\in\mathbb Z^4$, the resulting partition is $\mathbb Z^4$-periodic, with unit cell $\mathrm Y^4=[0,1)^4$.

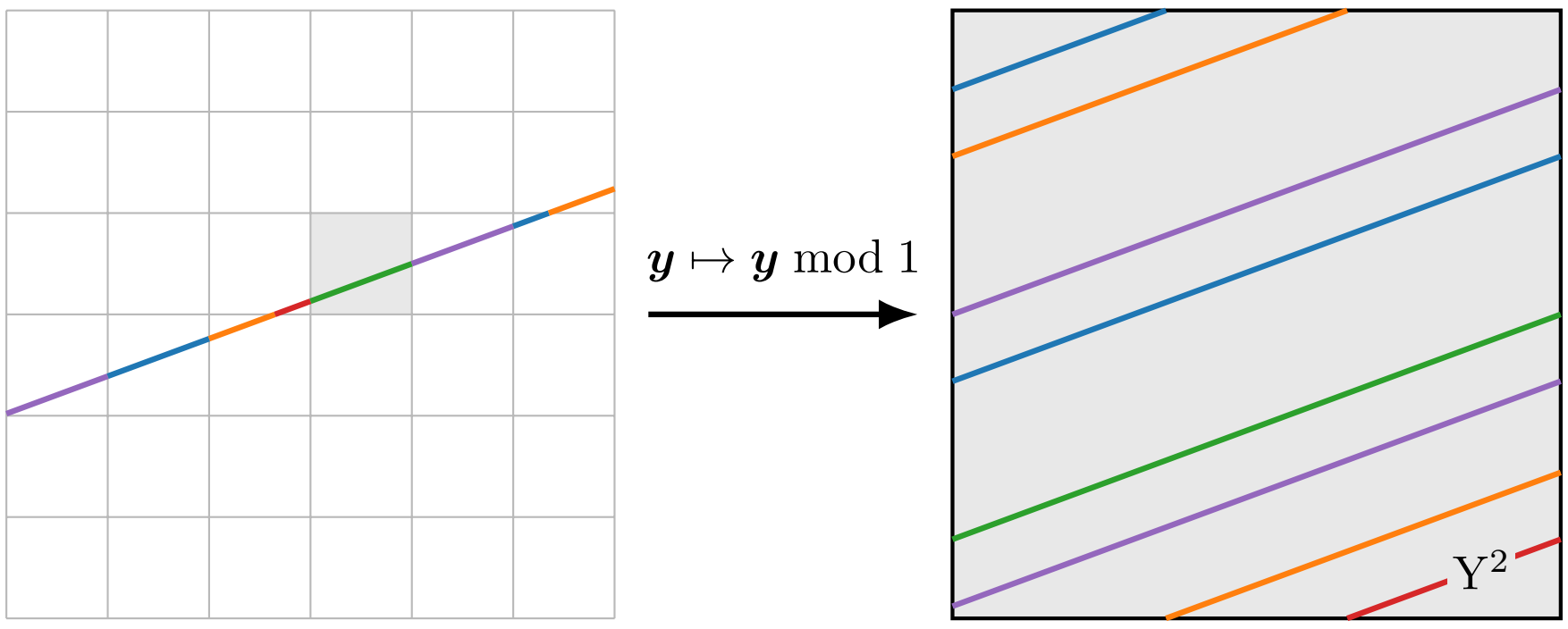
\begin{figure}[t!]
    \centering
    \begin{tikzpicture}[>=Latex]
        \begin{scope}[xshift=-3.5cm,x=0.75cm,y=0.75cm]
            \fill[gray!18] (0,0) rectangle (1,1);
            \draw[gray!55] (-3,-3) grid[step=1] (3,3);
            \draw[C4,very thick] (-3,-0.98)--(-2,-0.61);
            \draw[C0,very thick] (-2,-0.61)--(-1,-0.24);
            \draw[C1,very thick] (-1,-0.24)--(-0.351,0);
            \draw[C3,very thick] (-0.351,0)--(0,0.13);
            \draw[C2,very thick] (0,0.13)--(1,0.50);
            \draw[C4,very thick] (1,0.50)--(2,0.87);
            \draw[C0,very thick] (2,0.87)--(2.351,1);
            \draw[C1,very thick] (2.351,1)--(3,1.24);
        \end{scope}
        \node at (-3.5,-2.65) {(a) periodic square grid};
        \draw[->,very thick] (-1.0,0)--(1.0,0) node[midway,above=3pt] {$\bm y\mapsto\bm y\bmod 1$};
        \begin{scope}[xshift=1.25cm,yshift=-2.25cm,x=4.5cm,y=4.5cm]
            \fill[gray!18] (0,0) rectangle (1,1);
            \draw[thick] (0,0) rectangle (1,1);
            \draw[C4,very thick] (0,0.02)--(1,0.39);
            \draw[C0,very thick] (0,0.39)--(1,0.76);
            \draw[C1,very thick] (0,0.76)--(0.649,1);
            \draw[C3,very thick] (0.649,0)--(1,0.13);
            \draw[C2,very thick] (0,0.13)--(1,0.50);
            \draw[C4,very thick] (0,0.50)--(1,0.87);
            \draw[C0,very thick] (0,0.87)--(0.351,1);
            \draw[C1,very thick] (0.351,0)--(1,0.24);
            \node[fill=gray!18,inner sep=1pt] at (0.87,0.08) {$\mathrm Y^2$};
        \end{scope}
        \node at (3.5,-2.65) {(b) reference square};
    \end{tikzpicture}
    \caption{Two-dimensional analogue of the unit-cell mapping~\eqref{eq:pentagrid_plane_in_unit_cell}. (a) A line divided into colored segments by a periodic square grid. (b) The same segments translated by integer vectors into the reference square $\mathrm Y^2$.}
    \label{fig:irrational_cut}
\end{figure}

To identify the four-dimensional coordinates of $\bm q(\bm z)$, we introduce the linear map $\bm A^4:\mathbb R^2\to\mathbb R^4$, represented by the first four rows of $\bm A^5$, i.e., $\bm A^4=(\bm a_1, \bm a_2, \bm a_3, \bm a_4)^\mathsf{T}\in\mathbb R^{4\times2}$. Noting that $\bm A^5=\bm B\bm A^4$, substituting $\bm\gamma = \bm B\bm y_{\bm\gamma}$, with $\bm y_{\bm\gamma} = (\gamma_1,\gamma_2,\gamma_3,\gamma_4)^\mathsf{T}$, into equation~\eqref{eq:five_dimensional_point} yields
\begin{equation}
    \label{eq:plane_in_four_coordinates}
    \bm q(\bm z)
    =
    \bm B
    \left[
        \bm y_{\bm\gamma}
        +\bm A^4\bm z
    \right].
\end{equation}
Thus, $\bm y_{\bm\gamma}+\bm A^4\bm z$ is the four-dimensional coordinate vector of $\bm q(\bm z)$. Moreover, $\Pi_{\bm\gamma}^5=\bm B\Pi_{\bm\gamma}^4$, with the shifted plane represented in four-dimensional coordinates as
\begin{equation}
    \label{eq:pi_plane_4}
    \Pi_{\bm\gamma}^4
    =
    \bm y_{\bm\gamma}+\bm A^4\mathbb R^2.
\end{equation}

Due to the $\mathbb Z^4$-periodicity of $\{S_{\bm\kappa}^4\}_{\bm\kappa\in\mathbb Z^5}$, any point $\bm y\in\Pi_{\bm\gamma}^4$ can be mapped into the four-dimensional unit cell through
\begin{equation}
    \label{eq:pentagrid_plane_in_unit_cell}
    \bm y(\bm z)
    =
    \left(
        \bm y_{\bm\gamma}
        +\bm A^4\bm z
    \right)
    \bmod 1
    \in
    \mathrm Y^4.
\end{equation}
For $\bm z\in M$, equations~\eqref{eq:four_dimensional_partition} and \eqref{eq:plane_in_four_coordinates} show that $\bm y_{\bm\gamma} +\bm A^4\bm z\in S_{\bm\kappa(M)}^4$. Hence, the intersections $\Pi_{\bm\gamma}^4\cap S_{\bm\kappa(M)}^4$  correspond to the pentagrid elements and generate the same Penrose tiling as the dual construction in section~\ref{sec:crossing_to_rhombus}. Figure~\ref{fig:irrational_cut} illustrates a two-dimensional analogue of the unit-cell mapping~\eqref{eq:pentagrid_plane_in_unit_cell}.

\subsubsection{Distribution of material properties}

The four-dimensional representation developed above is now used to encode the distribution of elastic properties in a Penrose tiling plane. We assign the elasticity tensors $\bm{\mathsf C}_{\rm thick}$ and $\bm{\mathsf C}_{\rm thin}$ to the thick and thin rhombi, respectively. The resulting elasticity tensor in the physical plane is given by
\begin{equation}\label{eq:physical_material}
    \bm{\mathsf C}(\bm x)=\begin{cases}
    \bm{\mathsf C}_{\rm thick}, & \text{if $\bm x$ lies in a thick rhombus},\\
    \bm{\mathsf C}_{\rm thin}, & \text{if $\bm x$ lies in a thin rhombus}.
    \end{cases}
\end{equation}
To construct the corresponding periodic tensor on $\mathrm Y^4$, recall that the offset $\bm\gamma$ determines the shifted plane $\Pi_{\bm\gamma}^4$ and hence an entire Penrose tiling. Owing to the $\mathbb Z^4$-periodicity of the partition, different Penrose tilings may be sampled over a single unit cell by varying the offset from $\bm y_{\bm\gamma}$ to $\bm y_{\bm\gamma}+\bm y$, with $\bm y\in\mathrm Y^4$.

We now assign to each $\bm y$ the material property occupying a fixed point in the corresponding tiling, chosen as $\bm x=\bm0$. Since varying $\bm y$ changes the position of the tiling in the physical plane, we first remove the offset-dependent translation. Specifically, the tiling associated with the offset $\bm\gamma+\bm B\bm y$ has centered vertices
\begin{equation}\label{eq:centered_vertex}
    \bm p^\circ(M;\bm y)
    =
    \bm p(M)
    -
    (\bm A^5)^\mathsf{T}
    \big(
        \bm\gamma+\bm B\bm y
    \big)
    =
    (\bm A^5)^\mathsf{T}
    \big(
        \bm\kappa(M)-\bm\gamma-\bm B\bm y
    \big).
\end{equation}
The periodic elasticity tensor is then defined by the rhombus containing $\bm x=\bm0$ in this centered tiling,
\begin{equation}
    \label{eq:periodic_material}
    \bm{\mathsf C}_{\rm per}(\bm y)
    =
    \begin{cases}
        \bm{\mathsf C}_{\rm thick}
        & \text{if $\bm x=\bm0$ lies in a thick rhombus of the centered tiling},\\
        \bm{\mathsf C}_{\rm thin}
        & \text{if $\bm x=\bm0$ lies in a thin rhombus of the centered tiling},
    \end{cases}
\end{equation}
where we omit the dependence on the fixed offset $\bm \gamma$. Note that $\bm{\mathsf C}_{\rm per}$ is $\mathbb Z^4$-periodic,
since integer shifts of $\bm y$ induce identical shifts in the offset and pentagrid labels, leaving the centered vertices~\eqref{eq:centered_vertex} unchanged.

It remains to relate $\bm{\mathsf C}_{\rm per}$ to the prescribed physical distribution $\bm{\mathsf C}(\bm x)$, defined on the centered tiling with vertices~\eqref{eq:centered_vertex} at $\bm y=\bm0$. Consider offset deviations $\bm y=\bm A^4\bm z$, so that $\bm y_{\bm\gamma}+\bm y$ remains in $\Pi_{\bm\gamma}^4$ and the physical tiling remains unchanged apart from the centering translation. Substituting into~\eqref{eq:centered_vertex} gives
\begin{equation}
    \label{eq:centered_vertex_along_plane}
    \bm p^\circ
    \big(M;\bm A^4\bm z \big) =
    \bm p(M) -
    (\bm A^5)^\mathsf{T}
    \big(
        \bm\gamma+\bm B\bm A^4\bm z
    \big) =
    \bm p^\circ(M;\bm0)
    -
    (\bm A^5)^\mathsf{T}\bm A^5\bm z,
\end{equation}
where we have used $\bm A^5=\bm B\bm A^4$. Thus, the physical origin in the centered tiling associated with $\bm y=\bm A^4\bm z$ corresponds to the point $(\bm A^5)^\mathsf{T}\bm A^5\bm z$ in the centered tiling with offset $\bm\gamma$. Therefore, equation~\eqref{eq:centered_vertex} gives
\begin{equation}\label{eq:material_along_plane}
    \bm{\mathsf C}_{\rm per}(\bm A^4\bm z)=\bm{\mathsf C}\left((\bm A^5)^\mathsf{T}\bm A^5\bm z\right).
\end{equation}
Since $(\bm A^5)^\mathsf{T}\bm A^5=(5/2)\bm I_2$, we define the linear map $\RR_{\rm Pen}:\mathbb R^2\to\mathbb R^4$ by
\begin{equation}\label{eq:R_map}
    \RR_{\rm Pen}=\frac{2}{5}\bm A^4\in\mathbb R^{4\times2},
\end{equation}
Equation~\eqref{eq:material_along_plane} then yields
\begin{equation}\label{eq:quasiperiodic_material}
    \bm{\mathsf C}(\bm x)=\bm{\mathsf C}_{\rm per}(\RR_{\rm Pen}\,\bm x),
\end{equation}
a four-dimensional periodic form suitable for two-dimensional quasiperiodic homogenization~\citep{Wellander2019}.

\subsection{Three-dimensional icosahedral tiling}

We now consider the three-dimensional Ammann--Kramer--Neri tiling~\cite{Buganski2020}. By analogy with the pentagrid construction above, we consider a hexagrid consisting of six families of parallel planes, with prolate and oblate rhombohedra playing the role of the Penrose rhombi. The polyhedra formed by these planes receive six integer labels, which yield a periodic representation in $\mathbb R^6$~\cite{KramerNeri1984,DuneauKatz1985,Guenneau2020}.

\subsubsection{Three-dimensional hexagrid construction}

Let $\tau=(1+\sqrt{5})/2$ denote the golden ratio. We consider six unit vectors associated with the icosahedral directions, assembled into the linear map $\bm A^6:\mathbb R^3\to\mathbb R^6$, represented by ~\citep{Guenneau2020}
\begin{equation}\label{eq:icosahedral_directions}
    (\bm A^6)^\mathsf T=(\bm a_1\ \ \bm a_2\ \ \bm a_3\ \ \bm a_4\ \ \bm a_5\ \ \bm a_6)=\frac{1}{\sqrt{1+\tau^2}}\begin{pmatrix}\tau&\tau&0&-1&0&1\\0&0&1&\tau&1&\tau\\1&-1&-\tau&0&\tau&0\end{pmatrix}.
\end{equation}
The twelve points $\{\pm\bm a_j\}_{j=1}^6$ form the vertices of a regular icosahedron. These vectors satisfy
\begin{equation}\label{eq:icosahedral_properties}
    \|\bm a_j\|=1,\qquad |\bm a_i\cdot\bm a_j|=\frac{1}{\sqrt5},\quad i\ne j.
\end{equation}

We introduce an offset vector $\bm\gamma\in\mathbb R^6$ and, for $\bm z\in\mathbb R^3$, define the shifted projections
\begin{equation}\label{eq:icosahedral_projection}
    \bm q(\bm z)=\bm\gamma+\bm A^6\bm z.
\end{equation}
For every $k\in\mathbb Z$, the $k$th plane of family $j$ is given by
\begin{equation}\label{eq:icosahedral_planes}
    G_j^k=\{\bm z\in\mathbb R^3 \ : \ q_j(\bm z)=k\}.
\end{equation}
For fixed $j$, the collection $\{G_j^k\}_{k\in\mathbb Z}$ represents a family of parallel planes separated by unit distance. The union of the six families is called the \emph{icosahedral hexagrid}. We assume that the hexagrid is regular, meaning that no point lies on four or more planes.

\subsubsection{Dual tiling}

The hexagrid planes partition $\mathbb R^3$ into polyhedral regions, each referred to as a \emph{hexagrid element} $M$. Each element is assigned an integer label
\begin{equation}\label{eq:icosahedral_element_label}
\bm\kappa(M)=\big(\lceil q_1(\bm z)\rceil,\ldots,\lceil q_6(\bm z)\rceil\big)^\mathsf T\in\mathbb Z^6,\qquad \bm z\in M.
\end{equation}
The dual tiling is constructed by mapping each hexagrid element to a physical vertex $\bm p\in\mathbb R^3$,
\begin{equation}\label{eq:icosahedral_dual_vertex}
    \bm p(M)=\sum_{j=1}^6\kappa_j(M)\bm a_j=(\bm A^6)^\mathsf T\bm\kappa(M).
\end{equation}
Crossing a plane of family $j$ changes the integer label by $\pm\bm g_j$ and the corresponding physical vertex by $\pm\bm a_j$, where $\bm g_j$ denotes the $j$th standard basis vector of $\mathbb Z^6$.

\begin{figure}[t!]
    \centering
    \begin{minipage}[c]{0.46\textwidth}
        \centering
        \includegraphics[width=0.92\linewidth]{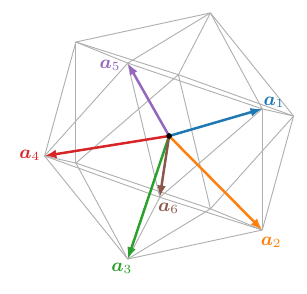}
    \end{minipage}
    \hspace{1mm}
    \begin{minipage}[c]{0.46\textwidth}
        \centering
        \includegraphics[width=0.82\linewidth]{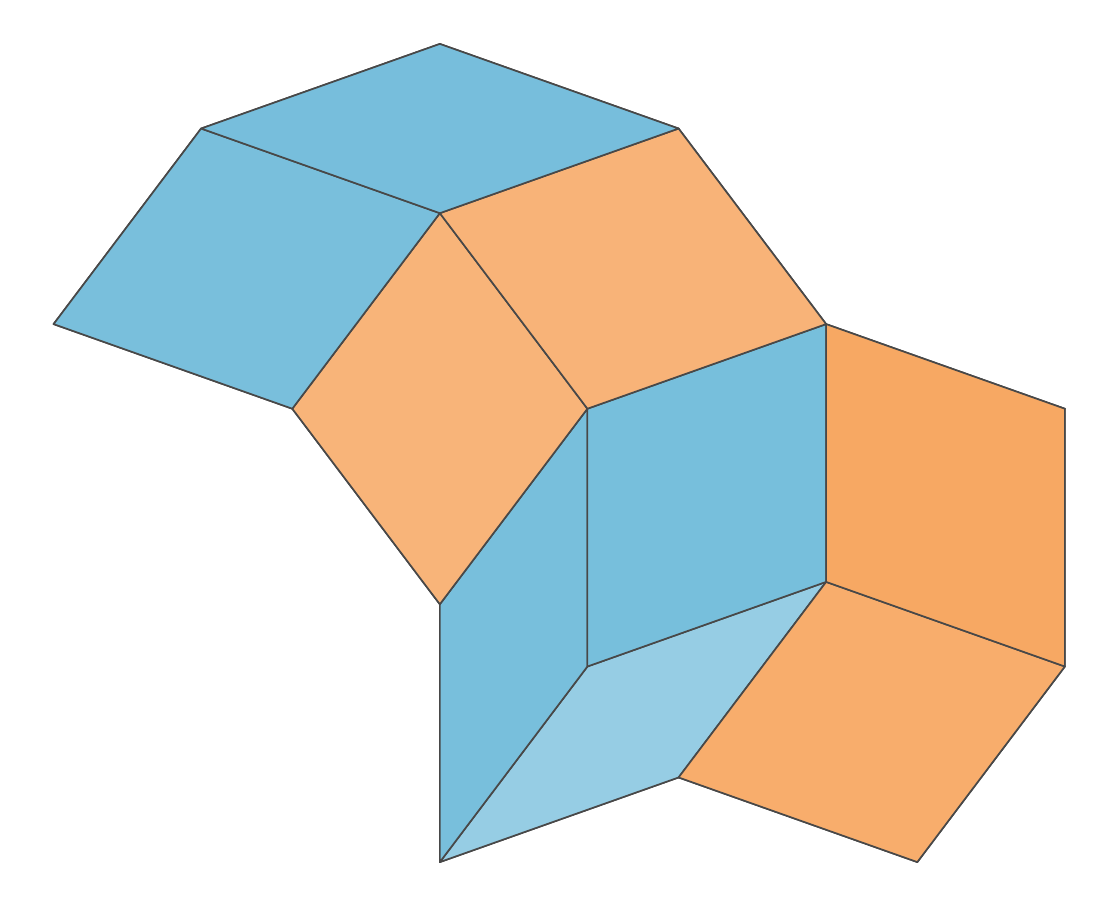}
    \end{minipage}
    \caption{Icosahedral directions and Ammann--Kramer--Neri tiling. The six vectors $\{\bm a_j\}_{j=1}^6$ point from the center of a regular icosahedron toward one vertex for each opposite pair (left). The finite patch consists of connected prolate and oblate rhombohedra, shown in blue and orange, respectively (right).}
    \label{fig:icosahedral_geometry}
\end{figure}

Suppose that planes from families $r$, $s$, and $t$ meet at a point. The vectors $\bm a_r$, $\bm a_s$, and $\bm a_t$ are linearly independent, so exactly eight hexagrid elements meet at the crossing. Denoting one of the corresponding physical vertices by $\bm p$, all eight vertices are given by
\begin{equation}\label{eq:icosahedral_rhombohedron}
    \left\{\bm p+d_r\bm a_r+d_s\bm a_s+d_t\bm a_t \ : \  d_r,d_s,d_t\in\{0,1\}\right\}.
\end{equation}
These vertices form a golden rhombohedron with unit edges parallel to $\bm a_r$, $\bm a_s$, and $\bm a_t$. The rhombohedron is either prolate or oblate, with volume ratio $\tau$ between the two types. Applying this construction to all hexagrid intersections generates the Ammann--Kramer--Neri tiling of the physical three-dimensional space~\citep{KramerNeri1984,DuneauKatz1985}. Figure~\ref{fig:icosahedral_geometry} shows the six icosahedral directions and a finite patch of this tiling.

\subsubsection{Six-dimensional periodic representation}

As $\bm z$ varies over $\mathbb R^3$, the points $\bm q(\bm z)\in\mathbb R^6$ given by~\eqref{eq:icosahedral_projection} form the shifted three-dimensional plane
\begin{equation}\label{eq:icosahedral_shifted_plane}
    \Pi_{\bm\gamma}^6=\{\bm q(\bm z) \ : \ \bm z\in\mathbb R^3\}=\bm\gamma+\bm A^6\mathbb R^3.
\end{equation}
For each $\bm\kappa\in\mathbb Z^6$, we identify an open unit hypercube
\begin{equation}\label{eq:icosahedral_hypercube}
    Q^6_{\bm\kappa}=\{\bm q\in\mathbb R^6 \ : \ \kappa_j-1<q_j<\kappa_j,\ j=1,\ldots,6\}.
\end{equation}
The hypercube closures form the standard $\mathbb Z^6$-periodic grid in $\mathbb R^6$. The shifted plane $\Pi_{\bm\gamma}^6$ passes through this periodic grid, and each hexagrid element $M$ is identified uniquely with a three-dimensional section,
\begin{equation}\label{eq:icosahedral_section}
    \{\bm q(\bm z) \ : \ \bm z\in M\}=\Pi_{\bm\gamma}^6\cap Q^6_{\bm\kappa(M)}.
\end{equation}
The six icosahedral directions~\eqref{eq:icosahedral_directions} are such that $(\bm A^6)^\mathsf{T}\bm n=\bm0$ for $\bm n\in\mathbb Z^6$ only when $\bm n=\bm0$. Thus, the six integer labels cannot be reduced to fewer independent coordinates as in the Penrose tiling.  Accordingly, we retain the six-dimensional periodic representation, with unit cell $\mathrm Y^6=[0,1)^6$.

Due to the $\mathbb Z^6$-periodicity of the hypercube grid, the shifted plane can be mapped into the unit cell through
\begin{equation}\label{eq:icosahedral_plane_in_unit_cell}
    \bm y(\bm z)=(\bm\gamma+\bm A^6\bm z)\bmod 1\in\mathrm Y^6.
\end{equation}
For $\bm z\in M$, equation~\eqref{eq:icosahedral_plane_in_unit_cell} translates $\bm q(\bm z)$ by an integer vector into $\mathrm Y^6$, while $\bm\kappa(M)$ in~\eqref{eq:icosahedral_element_label} identifies the hypercube $Q^6_{\bm\kappa(M)}$ containing $\bm q(\bm z)$. Equation~\eqref{eq:icosahedral_dual_vertex} finally maps this label to the physical vertex $\bm p(M)$.

\subsubsection{Distribution of material properties}

We now use the six-dimensional representation to encode the distribution of elastic properties in the icosahedral tiling. We assign the elasticity tensors ${\bm{\mathsf C}}_{\rm prolate}$ and ${\bm{\mathsf C}}_{\rm oblate}$ to the prolate and oblate rhombohedra, respectively. The resulting elasticity tensor in the physical space is written as
\begin{equation}\label{eq:icosahedral_physical_material}
    {\bm{\mathsf C}}(\bm x)=\begin{cases}{\bm{\mathsf C}}_{\rm prolate},&\text{if $\bm x$ lies in a prolate rhombohedron},\\{\bm{\mathsf C}}_{\rm oblate},&\text{if $\bm x$ lies in an oblate rhombohedron}.\end{cases}
\end{equation}

To construct the corresponding periodic tensor on $\mathrm Y^6$, recall that the offset $\bm\gamma$ determines the shifted plane $\Pi_{\bm\gamma}^6$ and hence an entire icosahedral tiling. Owing to the $\mathbb Z^6$-periodicity of the hypercube grid, different tilings may be sampled over a single unit cell by varying $\bm y\in\mathrm Y^6$. As in the Penrose construction, we assign to each $\bm y$ the material property occupying a fixed point in the corresponding tiling, chosen as $\bm x=\bm0$. Since varying the offset also changes the position of the tiling in the physical space, we first remove the offset-dependent translation. Specifically, the tiling associated with the offset $\bm\gamma+\bm y$ has centered vertices
\begin{equation}
    \label{eq:icosahedral_centered_vertex}
    \bm p^\circ(M;\bm y)
    =
    \bm p(M)
    -
    (\bm A^6)^\mathsf{T}
    \big(
        \bm\gamma+\bm y
    \big)
    =
    (\bm A^6)^\mathsf{T}
    \big(
        \bm\kappa(M)-\bm\gamma-\bm y
    \big).
\end{equation}
The periodic elasticity tensor is then defined by the rhombohedron containing $\bm x=\bm0$ in this centered tiling,
\begin{equation}\label{eq:icosahedral_periodic_material}
    {\bm{\mathsf C}}_{\rm per}(\bm y)=\begin{cases}{\bm{\mathsf C}}_{\rm prolate},&\text{if $\bm x=\bm0$ lies in a prolate rhombohedron},\\{\bm{\mathsf C}}_{\rm oblate},&\text{if $\bm x=\bm0$ lies in an oblate rhombohedron}.\end{cases}
\end{equation}
Note that $\bm{\mathsf C}_{\rm per}$ is $\mathbb Z^6$-periodic, since integer shifts of $\bm y$ induce identical shifts in the offset and hexagrid labels, leaving the centered vertices~\eqref{eq:icosahedral_centered_vertex} unchanged.

It remains to relate ${\bm{\mathsf C}}_{\rm per}$ to the prescribed physical distribution ${\bm{\mathsf C}}(\bm x)$, defined on the centered tiling with vertices~\eqref{eq:icosahedral_centered_vertex} at $\bm y=\bm0$. Consider offset deviations $\bm y=\bm A^6\bm z$ along the shifted plane $\Pi_{\bm\gamma}^6$. Substituting into equation~\eqref{eq:icosahedral_centered_vertex} gives
\begin{equation}
    \label{eq:icosahedral_centered_vertex_along_plane}
    \bm p^\circ
    \big(M;\bm A^6\bm z\big)
    =
    \bm p(M)
    -
    (\bm A^6)^\mathsf{T}
    \big(
        \bm\gamma+\bm A^6\bm z
    \big)
    =
    \bm p^\circ(M;\bm0)
    -
    (\bm A^6)^\mathsf{T}\bm A^6\bm z,
\end{equation}
Thus, the physical origin in the centered tiling associated with $\bm y=\bm A^6\bm z$ corresponds to the point $(\bm A^6)^\mathsf T\bm A^6\bm z$ in the centered tiling with offset $\bm\gamma$. Therefore,
\begin{equation}\label{eq:icosahedral_material_along_plane}
    {\bm{\mathsf C}}_{\rm per}(\bm A^6\bm z)={\bm{\mathsf C}}\left((\bm A^6)^\mathsf T\bm A^6\bm z\right).
\end{equation}
Since the six unit directions satisfy $(\bm A^6)^\mathsf T\bm A^6=2\bm I_3$, we define the linear map $\RR_{\rm AKN}:\mathbb R^3\to\mathbb R^6$ by
\begin{equation}\label{eq:icosahedral_R_map}
    \RR_{\rm AKN}=\frac12\bm A^6\in\mathbb R^{6\times3},
\end{equation}
Equation~\eqref{eq:icosahedral_material_along_plane} then yields
\begin{equation}\label{eq:icosahedral_quasiperiodic_material}
    {\bm{\mathsf C}}(\bm x)={\bm{\mathsf C}}_{\rm per}(\RR_{\rm AKN}\,\bm x).
\end{equation}
a six-dimensional periodic form suitable for three-dimensional quasiperiodic homogenization~\citep{Guenneau2020}.

\section{Homogenization}
\label{sec:hom}

\subsection{Problem setting}

Let $\Omega\subset\mathbb R^n$ denote the physical domain, with boundary $\partial\Omega=\Gamma_{\rm D}\cup\Gamma_{\rm N}$. We consider a sequence of increasingly fine quasiperiodic composites in $\Omega$, with $\epsilon$ denoting the tile edge length relative to a characteristic macroscopic length, and seek the continuum limit as $\epsilon\to0$.

For Penrose or AKN tilings, equations~\eqref{eq:quasiperiodic_material} and \eqref{eq:icosahedral_quasiperiodic_material} take the common form
\begin{equation}
    \label{eq:general_scaled_material}
    \bm{\mathsf C}_\epsilon(\bm x)
    =
    \bm{\mathsf C}_{\rm per}
    \bigg(
        \frac{\RR\bm x}{\epsilon}
    \bigg),
    \qquad
    \bm x\in\Omega,
\end{equation}
where $\bm{\mathsf C}_{\rm per}$ is periodic on $\mathrm Y^m=[0,1)^m$. Here, $\RR=\RR_{\rm Pen}$ for the Penrose material, with $(n,m)=(2,4)$, and $\RR=\RR_{\rm AKN}$ for the icosahedral material, with $(n,m)=(3,6)$. Both maps satisfy
\begin{equation}
    \label{eq:R_condition}
    \RR^\mathsf{T}\bm k\neq\bm0
    \qquad
    \forall\,\bm k\in\mathbb Z^m\setminus\{\bm0\}.
\end{equation}
We assume that $\bm{\mathsf C}_{\rm per}$ has the usual minor and major symmetries and that there exist constants $c_1,c_2>0$ such that
\begin{equation}
    \label{eq:elasticity_bounds}
    c_1|\bm\varepsilon|^2
    \leq
    \bm{\mathsf C}_{\rm per}(\bm y)\,\bm\varepsilon:\bm\varepsilon
    \leq
    c_2|\bm\varepsilon|^2
\end{equation}
for every $\bm\varepsilon\in\mathbb R_{\rm sym}^{n\times n}$ and almost every $\bm y\in\mathrm Y^m$. Letting $\bm\varepsilon(\bm u)=\operatorname{sym} \nabla\bm u$ denote the small-strain tensor, the corresponding $\epsilon$-scaled microscopic strain-energy density reads
\begin{equation}
    \label{eq:microscopic_density}
    W_\epsilon(\bm x,\bm\varepsilon) = W\bigg(
        \frac{\RR\bm x}{\epsilon},\bm\varepsilon\bigg)
    =
    \frac12
    \bm{\mathsf C}_{\rm per} \bigg(
        \frac{\RR\bm x}{\epsilon}
    \bigg)\bm\varepsilon
    :
    \bm\varepsilon,
\end{equation}
which yields the microscopic stress tensor
\begin{equation}
    \label{eq:microscopic_stress}
    \bm\sigma_\epsilon(\bm x,\bm\varepsilon)
    =
    \frac{\partial W_\epsilon}{\partial\bm\varepsilon}
    (\bm x,\bm\varepsilon)
    =
    \bm{\mathsf C}_{\rm per}\bigg(
        \frac{\RR\bm x}{\epsilon}
    \bigg)\bm\varepsilon.
\end{equation}

For a prescribed displacement $\bar{\bm u}$ on $\Gamma_{\rm D}$, the displacement field $\bm u_\epsilon$ minimizes
\begin{equation}
    \label{eq:microscopic_energy}
    \mathcal F_\epsilon(\bm u)
    =
    \begin{dcases}
        \displaystyle
        \int_\Omega
            W
            \bigg(
                \frac{\RR\bm x}{\epsilon},
                \bm\varepsilon(\bm u)
            \bigg)
        \,\mathrm d\bm x
        -
        \mathcal L(\bm u),
        & \bm u\in\mathrm H_{\bar{\bm u}}^1(\Omega;\mathbb R^n),\\
        +\infty,
        & \text{otherwise},
    \end{dcases}
\end{equation}
where the external loading functional reads
\begin{equation}
    \mathcal L(\bm u)
    =
    \int_\Omega \bm b\cdot\bm u\,\mathrm d\bm x
    +
    \int_{\Gamma_{\rm N}}
        \bar{\bm t}\cdot\bm u\,\mathrm dS.
\end{equation}
The corresponding Euler--Lagrange equations are
\begin{equation}
    \label{eq:microscopic_equilibrium}
    \begin{cases}
        -\operatorname{div}
        \bm\sigma_\epsilon(\bm x,\bm\varepsilon(\bm u_\epsilon))
        =
        \bm b
        & \text{in }\Omega,\\
        \bm u_\epsilon=\bar{\bm u}
        & \text{on }\Gamma_{\rm D},\\
        \bm\sigma_\epsilon(\bm x,\bm\varepsilon(\bm u_\epsilon))\,\bm n
        =
        \bar{\bm t}
        & \text{on }\Gamma_{\rm N}.
    \end{cases}
\end{equation}
        
Under condition~\eqref{eq:R_condition}, it follows that $\bm u_\epsilon\rightharpoonup\bm u_0$ in $\mathrm H^1(\Omega;\mathbb R^n)$~\cite[Theorem~2]{Wellander2019}, where $\bm u_0$ minimizes the homogenized energy functional
\begin{equation}
    \label{eq:homogenized_energy}
    \mathcal F_0(\bm u)
    =
    \begin{dcases}
        \displaystyle
        \int_\Omega
            W_0
            \big(
                \bm\varepsilon(\bm u)
            \big)
        \,\mathrm d\bm x
        -
        \mathcal L(\bm u),
        & \bm u\in\mathrm H_{\bar{\bm u}}^1(\Omega;\mathbb R^n),\\
        +\infty,
        & \text{otherwise},
    \end{dcases}
\end{equation}
with a homogenized strain-energy density characterized by the effective elasticity tensor $\bm{\mathsf C}_0$,
\begin{equation}
    \label{eq:homogenized_density}
    W_0(\bm\varepsilon)
    =
    \frac12
    \bm{\mathsf C}_0\bm\varepsilon
    :
    \bm\varepsilon,
\end{equation}
yielding the homogenized stress
\begin{equation}
    \label{eq:homogenized_stress}
    \bm\sigma_0(\bm\varepsilon) = \frac{\partial W_0}{\partial\bm\varepsilon}(\bm\varepsilon) = \bm{\mathsf C}_0\bm\varepsilon.
\end{equation}
The corresponding homogenized Euler--Lagrange equations are
\begin{equation}
    \label{eq:homogenized_equilibrium}
    \begin{cases}
        -\operatorname{div}
        \bm\sigma_0(\bm\varepsilon(\bm u_0))
        =
        \bm b
        & \text{in }\Omega,\\
        \bm u_0=\bar{\bm u}
        & \text{on }\Gamma_{\rm D},\\
        \bm\sigma_0(\bm\varepsilon(\bm u_0))\,\bm n
        =
        \bar{\bm t}
        & \text{on }\Gamma_{\rm N}.
    \end{cases}
\end{equation}
        
\subsection{Higher-dimensional unit-cell problem}

It remains to determine the effective elasticity tensor $\bm{\mathsf C}_0$. For a constant macroscopic strain $\bm\varepsilon_0\in\mathbb R^{n\times n}_{\rm sym}$, prescribe an affine displacement $\bm u=\bm\varepsilon_0\bm x$ on $\partial\Omega$ and write the interior displacement as $\bm u=\bm\varepsilon_0\bm x+\bm v$, with $\bm v\in\mathrm H^1_0(\Omega;\mathbb R^n)$ denoting a fluctuation. The corresponding minimum microscopic energy is given by
\begin{equation}\label{eq:min_energy}
    m_\epsilon(\bm\varepsilon_0;\Omega)
    =
    \min_{\bm v\in\mathrm H^1_0(\Omega;\mathbb R^n)}
    \int_\Omega
    W\bigg(\frac{\bm Q\bm x}{\epsilon},\bm\varepsilon_0+\bm\varepsilon(\bm v)\bigg)
    \,\mathrm d\bm x.
\end{equation}
Convergence of the minimum energies gives, as in standard variational homogenization~\cite{Dalmaso1993},
\begin{equation}\label{eq:affine_minimum_limit}
    \lim_{\epsilon\to0}m_\epsilon(\bm\varepsilon_0;\Omega)
    =
    \min_{\bm v\in\mathrm H^1_0(\Omega;\mathbb R^n)}
    \int_\Omega W_0\big(\bm\varepsilon_0+\bm\varepsilon(\bm v)\big)\,\mathrm d\bm x.
\end{equation}
The minimum on the right-hand side is attained at $\bm v=\bm0$ and therefore equals $|\Omega|W_0(\bm\varepsilon_0)$. Consequently,
\begin{equation}\label{eq:affine_limit}
    W_0(\bm\varepsilon_0)
    =
    \frac1{|\Omega|}
    \lim_{\epsilon\to0}m_\epsilon(\bm\varepsilon_0;\Omega).
\end{equation}

To express the limit~\eqref{eq:affine_limit} on the periodic unit cell, consider a smooth $\mathrm Y^m$-periodic vector field $\bm w:\mathbb R^m\to\mathbb R^n$ with zero mean and the displacement fluctuation
\begin{equation}
    \bm v_\epsilon(\bm x)=\epsilon\bm w\bigg(\frac{\bm Q\bm x}{\epsilon}\bigg).
\end{equation}
The chain rule gives the projected symmetric gradient
\begin{equation}\label{eq:projected_strain}
    \bm\varepsilon(\bm v_\epsilon)(\bm x)=\bm\varepsilon_{\bm Q}(\bm w)\bigg(\frac{\bm Q\bm x}{\epsilon}\bigg),\qquad
    \bm\varepsilon_{\bm Q}(\bm w)=\operatorname{sym}\big((\nabla_{\bm y}\bm w)\bm Q\big).
\end{equation}
Thus, differentiation in the $m$-dimensional unit cell occurs only along the $n$ physical directions selected by $\bm Q$. Moreover, condition~\eqref{eq:R_condition} ensures that spatial averages along the physical section converge to the corresponding unit-cell averages. For every sufficiently regular $\mathrm Y^m$-periodic function $f$,
\begin{equation}\label{eq:cut_project_average}
    \lim_{\epsilon\to0}
    \frac{1}{|\Omega|}
    \int_\Omega
    f\bigg(\frac{\bm Q\bm x}{\epsilon}\bigg)
    \,\mathrm d\bm x
    =
    \frac{1}{|\mathrm Y^m|}
    \int_{\mathrm Y^m}f(\bm y)\,\mathrm d\bm y.
\end{equation}

Applying~\eqref{eq:cut_project_average} to $f(\bm y)=W\big(\bm y,\bm\varepsilon_0+\bm\varepsilon_{\bm Q}(\bm w)\big)$ gives the mean energy associated with the periodic fluctuation $\bm w$. The elasticity unit-cell formulation presented in~\citep[Theorem~2]{Wellander2019} admits the following quadratic variational characterization:
\begin{equation}\label{eq:cell_energy}
    W_0(\bm\varepsilon_0)
    =
    \inf_{\bm w\in\mathrm H^1_{\rm per}(\mathrm Y^m;\mathbb R^n)}
    \frac{1}{|\mathrm Y^m|}
    \int_{\mathrm Y^m}
    W\big(\bm y,\bm\varepsilon_0+\bm\varepsilon_{\bm Q}(\bm w)\big)
    \,\mathrm d\bm y.
\end{equation}
The associated stationarity condition gives the unit-cell problem in weak form,
\begin{equation}\label{eq:higher_dimensional_corrector_problem}
    \int_{\mathrm Y^m}
    \bm{\mathsf C}_{\rm per}(\bm y)
    \big(\bm\varepsilon_0+\bm\varepsilon_{\bm Q}(\bm w)\big)
    :
    \bm\varepsilon_{\bm Q}(\tilde{\bm w})
    \,\mathrm d\bm y
    =
    0
    \qquad
    \forall\,\tilde{\bm w}\in\mathrm H^1_{\rm per}(\mathrm Y^m;\mathbb R^n).
\end{equation}
Here, $\bm\varepsilon_{\bm Q}(\bm w)$ denotes the generalized projected strain corrector; the corrector itself need not belong to $\mathrm H^1_{\rm per}$~\cite[Proposition~14]{Wellander2019}. The corresponding equilibrium equations are formally written as
\begin{equation}\label{eq:higher_dimensional_corrector_strong}
    \begin{cases}
        -\operatorname{div}_{\bm y}
        \big\{
            \bm{\mathsf C}_{\rm per}(\bm y)
            [
                \bm\varepsilon_0+\bm\varepsilon_{\bm Q}(\bm w)
            ]
            \bm Q^\mathsf T
        \big\}
        =
        \bm0
        & \text{in }\mathrm Y^m,\\
        \bm w\text{ is periodic}
        & \text{on }\partial\mathrm Y^m,\\
        \big\{
            \bm{\mathsf C}_{\rm per}(\bm y)
            [
                \bm\varepsilon_0+\bm\varepsilon_{\bm Q}(\bm w)
            ]
            \bm Q^\mathsf T
        \big\}\bm n_{\bm y}
        \text{ is anti-periodic}
        & \text{on }\partial\mathrm Y^m,\\
        \displaystyle\int_{\mathrm Y^m}\bm w\,\mathrm d\bm y=\bm0.&
    \end{cases}
\end{equation}
The zero-mean condition fixes the arbitrary additive constant in the corrector.

To determine $\bm{\mathsf C}_0$, we consider the $n(n+1)/2$ independent macroscopic strain modes
\begin{equation}\label{eq:macroscopic_strain_modes}
    \bm\varepsilon_0^{(ij)}
    =
    \frac12
    \big(
        \bm e_i\otimes\bm e_j+
        \bm e_j\otimes\bm e_i
    \big),
    \qquad
    1\leq i\leq j\leq n,
\end{equation}
where $\{\bm e_i\}_{i=1}^n$ is the standard basis of $\mathbb R^n$. If $\bm\varepsilon_{\bm Q}(\bm w^{(ij)})$ denotes the projected strain solution of~\eqref{eq:higher_dimensional_corrector_problem} for $\bm\varepsilon_0=\bm\varepsilon_0^{(ij)}$, then
\begin{equation}\label{eq:effective_elasticity_components}
    (\bm{\mathsf C}_0)_{ijkl}
    =
    \frac{1}{|\mathrm Y^m|}
    \int_{\mathrm Y^m}
    \left[
        \bm{\mathsf C}_{\rm per}(\bm y)
        \big(
            \bm\varepsilon_0^{(ij)}+
            \bm\varepsilon_{\bm Q}(\bm w^{(ij)})
        \big)
    \right]_{kl}
    \,\mathrm d\bm y.
\end{equation}
Thus, three higher-dimensional cell problems determine the effective tensor of the Penrose material, while six determine that of the icosahedral material. These problems are periodic in $m$ dimensions but involve derivatives only along the $n$ physical directions encoded by $\bm Q$.

\section{Computational implementation}
\label{sec:comp_imp}

The evaluation of the homogenized elasticity tensor $\bm{\mathsf C}_0$ through~\eqref{eq:effective_elasticity_components} requires the solution of the unit-cell problem~\eqref{eq:higher_dimensional_corrector_strong} for each independent macroscopic strain mode~\eqref{eq:macroscopic_strain_modes}. We exploit the periodicity of this problem in $\mathrm Y^m$ through an FFT-based solution strategy~\citep{moulinec1998numerical,vondvrejc2014fft,lucarini2019dbfft}, here applied to the higher-dimensional problem. For the Penrose material, we introduce a truncation scheme for the FFT solver and a material sampling scheme for constructing the higher-dimensional unit cell, both designed to preserve the fivefold rotational symmetry. Then, we present the construction of the higher-dimensional unit cells and the numerical evaluation of $\bm{\mathsf C}_0$.


\subsection{FFT-based solution of the unit-cell problems}

\subsubsection{Fourier formulation}

The $\mathrm Y^m$-periodic corrector $\bm w$ can be expressed in Fourier form as
\begin{equation}\label{eq:fourier_corrector}
\bm w(\bm y)=\sum_{\bm k\in\mathbb Z^m}\widehat{\bm w}(\bm k)\exp(2\pi i\bm k\cdot\bm y).
\end{equation}
Here, $\bm k\in\mathbb Z^m$ is the higher-dimensional wavevector and $\widehat{\bm w}(\bm k)=\mathscr F[\bm w](\bm k)$ is the Fourier coefficient, where $\mathscr F$ denotes the Fourier transform. In view of~\eqref{eq:projected_strain}, the projected symmetric gradient $\bm\varepsilon_{\bm Q}$ in~\eqref{eq:higher_dimensional_corrector_strong} has the Fourier representation
\begin{equation}\label{eq:fft_strain}
    \mathscr F[{\bm\varepsilon_{\bm Q}}(\bm w)](\bm k)
    =\frac{i}{2}\, \big(\widehat{\bm w}(\bm k)\otimes\bm \xi(\bm k)+\bm \xi(\bm k)\otimes\widehat{\bm w}(\bm k)\big),
\end{equation}
where $\bm \xi(\bm k)=2\pi\RR^{\mathsf T}\bm k$ denotes the projected wavevector. The projected strain field is then recovered in the higher-dimensional unit cell by an inverse Fourier transform and multiplied pointwise by the elasticity tensor to obtain the corrector-induced stress,
\begin{equation}\label{eq:fft_stress}
    \bm\sigma_{\bm w}(\bm y)
    =\Cper(\bm y)\,\bm\varepsilon_{\bm Q}(\bm w)(\bm y).
\end{equation}
The Fourier transform of the projected divergence of the corrector-induced stress,  $\operatorname{div}_{\bm y}\{\bm{\mathsf C}_{\rm per}(\bm y)[\bm\varepsilon_{\bm Q}(\bm w)]\bm Q^\mathsf T\}$, takes the form
\begin{equation}\label{eq:fft_corrector_divergence}
    \mathscr F\!\left[\operatorname{div}_{\bm y}(\bm\sigma_{\bm w}\bm Q^{\mathsf T})\right](\bm k)
    =i\,\mathscr F[\bm\sigma_{\bm w}](\bm k)\,\bm \xi(\bm k).
\end{equation}
Likewise, for the stress induced by the prescribed macroscopic strain,
\begin{equation}\label{eq:fft_macro_divergence}
    \mathscr F\!\left[\operatorname{div}_{\bm y}\bigl((\Cper\bm\varepsilon_0)\bm Q^{\mathsf T}\bigr)\right](\bm k)
    =i\,\mathscr F[\Cper\bm\varepsilon_0](\bm k)\,\bm \xi(\bm k).
\end{equation}
Applying the inverse Fourier transform and substituting these expressions into the unit-cell problem~\eqref{eq:higher_dimensional_corrector_strong} gives the linear system
\begin{equation}\label{eq:fft_linear_system}
    \bm H\bm w=\bm b,
\end{equation}
where
\begin{equation}\label{eq:fft_operator}
    \bm H\bm w
    =\mathscr F^{-1}\!\left[-i\,\mathscr F[\bm\sigma_{\bm w}](\bm k)\,\bm \xi(\bm k)\right], \qquad 
    \bm b =\mathscr F^{-1}\!\left[i\,\mathscr F[\Cper\bm\varepsilon_0](\bm k)\,\bm \xi(\bm k)\right].
\end{equation}

The linear system~\eqref{eq:fft_linear_system} is solved iteratively without assembling $\bm H$. Each evaluation of $\bm H\bm w$ uses Fourier differentiation and pointwise multiplication by $\Cper$ in the higher-dimensional unit cell, with the forward and inverse transforms computed by multidimensional FFTs. The infinite Fourier series~\eqref{eq:fourier_corrector} must therefore be truncated to a finite set of modes, as described next.

\subsubsection{Fourier truncation and symmetry preservation}
\label{sec:fourier_truncation}

The FFT implementation requires truncating the Fourier series~\eqref{eq:fourier_corrector} to a finite set of modes. The retained set must be invariant under the Fourier-index transformations associated with the rotational symmetry of each material. As shown below, this requirement leads to a modified truncation criterion for the Penrose material and a conventional cutoff in the icosahedral case.

The spatial sampling index set is defined as
\begin{equation}
    \Lambda_N^m  = \{0,\ldots,N-1\}^m,
\end{equation}
where $N$ denotes the number of grid points along each coordinate direction. For odd $N$, the spatial sampling grid determines the corresponding FFT grid
\begin{equation}
    \hat{\Lambda}_N^m =
    \left\{
    \bm k \in \mathbb Z^m \ : \
    |k_\alpha| \leq \frac{N-1}{2},
    \ \alpha = 1,\ldots,m
    \right\}.
\end{equation}

\paragraph{Fourier truncation scheme for Penrose material} For a Fourier index $\bm k=(k_1,k_2,k_3,k_4)^\mathsf T\in\mathbb Z^4$, the corresponding wavevector in the physical plane follows from~\eqref{eq:five_directions} and~\eqref{eq:R_map} as
\begin{equation}
    \bm\xi(\bm k)
    =
    2\pi\bm Q^\mathsf T\bm k
    =
    \frac{4\pi}{5}
    \sum_{j=1}^{4}k_j\bm a_j.
\end{equation}
A physical rotation by $72^\circ$ maps each Penrose direction to the next one and therefore gives
\begin{equation}
    \bm R_{72}\,\bm\xi(\bm k)
    =
    \frac{4\pi}{5}
    \left(
        k_1\bm a_2
        +k_2\bm a_3
        +k_3\bm a_4
        +k_4\bm a_5
    \right), \qquad     \bm R_{72}
    =
    \begin{pmatrix}
        \cos 72^\circ & -\sin 72^\circ\\
        \sin 72^\circ & \cos 72^\circ
    \end{pmatrix}.
\end{equation}
Using the identity $\sum_{j=1}^{5}\bm a_j=\bm0$,
we have $\bm a_5=-(\bm a_1+\bm a_2+\bm a_3+\bm a_4)$. Substituting this relation yields
\begin{equation}
    \bm R_{72}\, \bm \xi(\bm k)
    =\frac{4\pi}{5}\big(
    -k_4\bm a_1
    +(k_1-k_4)\bm a_2
    +(k_2-k_4)\bm a_3
    +(k_3-k_4)\bm a_4
    \big).
\end{equation}
The physical rotation is therefore represented in Fourier-index space by
\begin{equation}\label{eq:transform}
    \bm k'
    =
    \begin{pmatrix}
        -k_4\\
        k_1-k_4\\
        k_2-k_4\\
        k_3-k_4
    \end{pmatrix}
    =
    \bm T\bm k,
    \qquad
    \bm T
    =
    \begin{pmatrix}
        0&0&0&-1\\
        1&0&0&-1\\
        0&1&0&-1\\
        0&0&1&-1
    \end{pmatrix},
\end{equation}
such that $\bm R_{72}\,\bm\xi(\bm k)=\bm\xi(\bm k')$. 

Although the physical rotation preserves the wavevector magnitude, $|\bm\xi(\bm k')|=|\bm\xi(\bm k)|$, it does not generally preserve the Euclidean norm of the corresponding Fourier index. Indeed, equation~\eqref{eq:transform} gives
\begin{equation}
    |\bm k'|^2-|\bm k|^2
    =
    3k_4^2-2k_4(k_1+k_2+k_3).
\end{equation}
Consequently, a conventional cutoff based on $|\bm k|$ may retain a Fourier mode while excluding its rotated counterpart, thereby breaking the fivefold symmetry of the truncated Fourier representation.

To preserve the fivefold symmetry of the Penrose construction, we introduce a truncation criterion based on a projected Fourier-index norm that is invariant under the index transformation associated with a $72^\circ$ rotation. To this end, we embed the four-dimensional Fourier index $\bm k$ in $\mathbb R^5$ as
\begin{equation}
    \bm k^{(5)}
    =
    (k_1,k_2,k_3,k_4,0)^\mathsf T,
\end{equation}
and project $\bm k^{(5)}$ onto the zero-sum hyperplane $\Xi^4$ defined in~\eqref{eq:Xi^4_definition},
\begin{equation}
    \bm k_\Xi=\bm P_\Xi\,\bm k^{(5)}, \qquad \bm P_\Xi
    =
    \bm I_5-\frac{1}{5}\bm 1_5\bm 1_5^\mathsf T,
\end{equation}
where $\bm 1_5=(1,1,1,1,1)^\mathsf T$. The resulting squared norm is
\begin{equation}\label{eq:projected_norm}
    |\bm k_\Xi|^2
    =
    \sum_{j=1}^{4}k_j^2
    -
    \frac{1}{5}
    \Bigg(
        \sum_{j=1}^{4}k_j
    \Bigg)^2.
\end{equation}

To verify that this norm is invariant under a
$72^\circ$ rotation, we see from~\eqref{eq:transform} that
\begin{equation}
    \sum_{j=1}^{4}k_j'
    =
    \sum_{j=1}^{4}k_j-5k_4, \qquad 
    |\bm k'|^2
    =
    |\bm k|^2+5k_4^2
    -2k_4\sum_{j=1}^{4}k_j.
\end{equation}
Substituting these expressions into~\eqref{eq:projected_norm} gives
\begin{equation}
    \begin{aligned}
        |\bm k'_\Xi|^2
        &=
        |\bm k'|^2
        -\frac{1}{5}
        \Bigg(
            \sum_{j=1}^{4}k_j'
        \Bigg)^2=
        |\bm k|^2
        -\frac{1}{5}
        \Bigg(
            \sum_{j=1}^{4}k_j
        \Bigg)^2
        =
        |\bm k_\Xi|^2.
    \end{aligned}
\end{equation}
Thus, the projected norm assigns the same value to
Fourier indices related by a $72^\circ$ rotation.

We retain the nonzero Fourier modes satisfying the invariant cutoff and require their complete rotational orbit to lie within the discrete Fourier grid. Accordingly,
\begin{equation}\label{eq:penrose_cutoff}
    \hat{\Lambda}_{\mathrm{Pen}}
    =
    \left\{
        \bm k\in\mathbb Z^4 \ : \
        0<|\bm k_\Xi|^2\leq\rho_{\mathrm{Pen}}^2,
        \
        \bm T^r\bm k\in \hat{\Lambda}_N^4,
        \
        r=0,\ldots,4
    \right\},
\end{equation}
where $\rho_{\mathrm{Pen}}>0$ is the prescribed cutoff radius that controls the number of retained modes. The orbit condition ensures that $\hat{\Lambda}_{\mathrm{Pen}}$ is closed under the index transformation~\eqref{eq:transform} and therefore preserves fivefold rotational symmetry.

\paragraph{Fourier truncation scheme for the icosahedral material} For the AKN tiling, the six Fourier coordinates correspond to the six pairs of opposite icosahedral directions. Icosahedral rotations act as signed permutations of these coordinates, preserving both their Euclidean norm and the finite FFT grid. Therefore, unlike the Penrose case, no additional symmetry-preserving truncation scheme is required. Nevertheless, we may also retain only the modes within a spherical cutoff, i.e.,
\begin{equation}\label{eq:akn_truncation}
    \hat{\Lambda}_{\mathrm{AKN}}
    =
    \left\{
    \bm k\in\mathbb Z^6 \ : \ 
    0<|\bm k|^2\leq\rho_{\mathrm{AKN}}^2,\
    \bm k\in \hat{\Lambda}_N^6
    \right\},
\end{equation}
where $\rho_{\mathrm{AKN}}>0$ is the prescribed cutoff radius.

\subsection{Higher-dimensional unit cell construction}

To implement the FFT solver, the periodic elasticity tensor $\bm{\mathsf C}_{\mathrm{per}}$ is sampled on a uniform grid in $\mathrm Y^m$, indexed by the spatial sampling set $\Lambda_N^m$. The reconstructed strain is evaluated on the spatial grid, allowing the stress~\eqref{eq:fft_stress} to be computed pointwise.

\paragraph{Unit cell construction for the Penrose material} 
The four-dimensional periodic unit
cell $\mathrm Y^4=[0,1)^4$ is sampled at
\begin{equation}
    \bm y_{\bm \alpha}
    = \frac{\bm \alpha + \frac15\bm 1_4}{N},
    \qquad \bm \alpha \in \Lambda_N^4,
\end{equation}
where $\bm 1_4 =(1,1,1,1)^\mathsf T$. In the five-dimensional representation, a $72^\circ$ rotation in physical space corresponds to a cyclic permutation of the five parent-cell coordinates. Under the zero-sum condition, this permutation induces the following transformation on the four-dimensional coordinates:
\begin{equation}
    \bm P\bm{y}
    = (y_2,y_3,y_4,-y_1-y_2-y_3-y_4)^{\mathsf T}
     \bmod 1.
\end{equation}
The choice of shift $1/5$ in the sampling grid ensures that $\bm P$ maps the grid onto itself. Since the thick- and thin-rhombus classification is invariant under a $72^\circ$ rotation, the five grid points $\bm y,\bm P\bm y,\ldots,\bm P^4\bm y$ are assigned the same material tensor, i.e.,
\begin{equation}\label{eq:penrose_orbit_invariance}
    \Cper(\bm P^r\bm y)
    =
    \Cper(\bm y),
    \qquad
    r=1,\ldots,4.
\end{equation}
Accordingly, the material classification is evaluated at one of these five points following section~\ref{sec:geometry} and assigned to the remaining four points.

\paragraph{Unit cell construction for the icosahedral material} 
For the AKN tiling, the full six-dimensional parent cell is sampled at
\begin{equation}
    \bm y_{\bm \alpha}
    = \frac{\bm \alpha + \frac12\bm 1_6}{N},
    \qquad \bm \alpha \in \Lambda_N^6,
\end{equation}
where $\bm 1_6=(1,1,1,1,1,1)^\mathsf T$. The material classification is evaluated at each grid point following section~\ref{sec:geometry} to construct the discrete elasticity field $\Cper$. The resulting $N^6$ material values provide the input for the FFT-based solution of the unit-cell problem for the AKN material.

\subsection{Homogenized elasticity tensor}

With both the FFT framework and the discrete elasticity field $\Cper$ established, the unit-cell problem in Fourier space~\eqref{eq:fft_linear_system} can be solved for $\bm w$ discretely. The homogenized elasticity tensor is then given by~\eqref{eq:effective_elasticity_components}, where the unit-cell integral is approximated by the arithmetic average of the corrected microscopic stress over the uniform spatial grid,
\begin{equation}
    (\bm{\mathsf{C}}_0)_{ijkl}
    \approx
    \frac{1}{N^m}
    \sum_{\bm \alpha\in\Lambda_N^m}
    \left[
    \bm{\mathsf{C}}_{\mathrm{per}}(\bm y_{\bm \alpha}):
    \left(
    \bm\varepsilon_0^{(ij)}
    +
    \bm\varepsilon_Q(\bm w^{(ij)})(\bm y_{\bm \alpha})
    \right)
    \right]_{kl}.
\end{equation}

\section{Numerical results}
\label{sec:examples}

This section presents numerical results for the proposed homogenization framework for quasiperiodic elastic materials. We compare the mechanical responses obtained from direct finite element simulations that
explicitly resolve the quasiperiodic microstructures
with those predicted by the corresponding homogenized
models. In particular, we investigate the convergence
of the strain energy toward the homogenized limit as
the characteristic microstructural scale $\epsilon$
decreases.

The construction of the higher-dimensional periodic unit cells and the FFT-based solution of the projected cell problems are implemented in \texttt{MATLAB}, while the direct and homogenized finite element models are created using the open-source \texttt{FEniCS} library in \texttt{Python}. Both the \texttt{MATLAB} and \texttt{Python} codes are provided as supplementary material.

\subsection{Two-dimensional Penrose tiling}

We compute the homogenized elasticity tensor for the Penrose material by solving the projected cell problems with the FFT-based method using $N=101$ grid points along each of the four Fourier directions. We assign Young's moduli $E^{(1)}=8$~MPa and $E^{(2)}=1$~MPa to the two rhombus types, with a common Poisson's ratio $\nu=0.30$. By way of illustration, when the stiff material is assigned to the thick rhombi and the soft material to the thin rhombi, we obtain the following homogenized elasticity tensor under plane-stress conditions in engineering notation:
\begin{equation}\label{eq:penrose_Chom}
    [\bm{\mathsf C}^{\rm eng}_0] =
    \begin{pmatrix}
    4.387 & 1.295 & 0 \\
    1.295 & 4.387 & 0 \\
    0 & 0 & 1.546
    \end{pmatrix}\,\mathrm{MPa}.
\end{equation}
The tensor exhibits an isotropic response in the plane, as expected from the fivefold symmetry.

\begin{figure}[t!]
    \vspace*{1em}
    \begin{center}
    \begin{tabular}{cc}
    \begin{overpic}[width=0.43\textwidth]{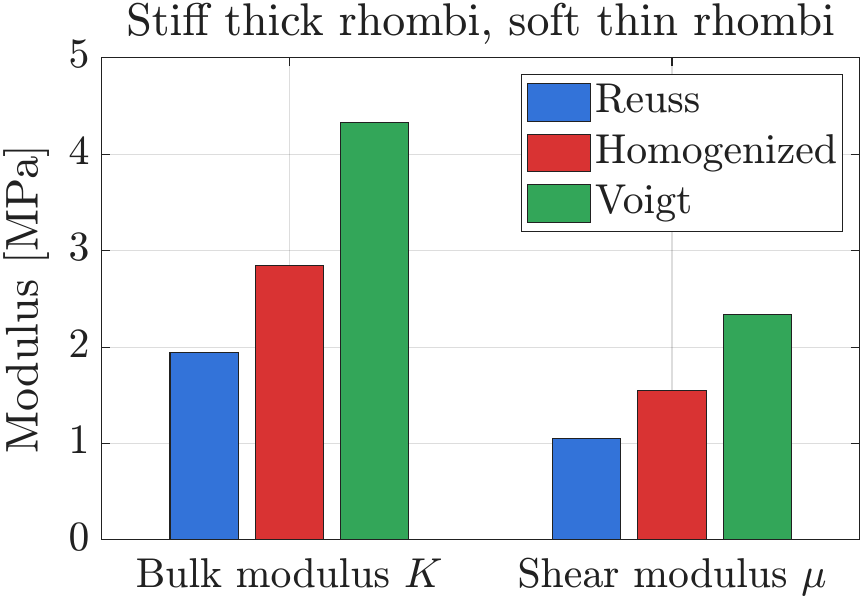}
        \put(0,2){\small{(a)}}
    \end{overpic} &
    \begin{overpic}[width=0.43\textwidth]{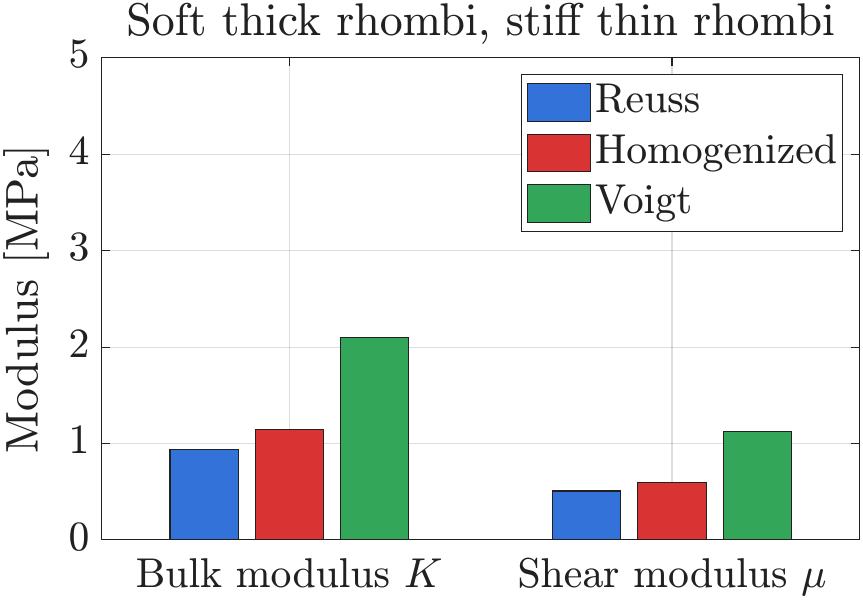}
        \put(0,2){\small{(b)}}
    \end{overpic}
    \end{tabular}
    \end{center}
    \vspace*{-1em}
    \caption{Comparison of the homogenized bulk and shear moduli
    with the Voigt and Reuss bounds for the Penrose tiling with different material assignments.
    (a) Stiff thick rhombi and soft thin rhombi.
    (b) Soft thick rhombi and stiff thin rhombi.}
    \label{fig:voigt_reuss_bounds}
\end{figure}

Figure~\ref{fig:voigt_reuss_bounds} compares the homogenized bulk and shear moduli of the Penrose tiling with the Voigt and Reuss bounds for two different material assignments. The bounds are evaluated using the area fractions of the thick ($f_{\mathrm{thick}}=\tau^2/(1+\tau^2) \approx 72.36\%$) and thin ($f_{\mathrm{thin}}=1/(1+\tau^2) \approx 27.64\%$) rhombi in the infinite tiling. When the thick rhombi are assigned the stiffer material, the homogenized bulk and shear moduli follow as $K_0=2.841$~MPa and $\mu_0=1.546$~MPa, respectively. These values lie between the corresponding Reuss bounds, $K_{\mathrm{R}}=1.947$~MPa and $\mu_{\mathrm{R}}=1.048$~MPa, and Voigt bounds, $K_{\mathrm{V}}=4.332$~MPa and $\mu_{\mathrm{V}}=2.333$~MPa. Inverting the phase assignment yields $K_0=1.147$~MPa and $\mu_0=0.598$~MPa, which similarly lie between the corresponding Reuss bounds, $K_{\mathrm{R}}=0.942$~MPa and $\mu_{\mathrm{R}}=0.507$~MPa, and Voigt bounds $K_{\mathrm{V}}=2.096$~MPa and $\mu_{\mathrm{V}}=1.129$~MPa. Thus, in both configurations, the homogenized moduli satisfy the Voigt--Reuss bounds, providing a first consistency check.

\begin{figure}[b!]
    \centering
    \includegraphics[width=\linewidth]{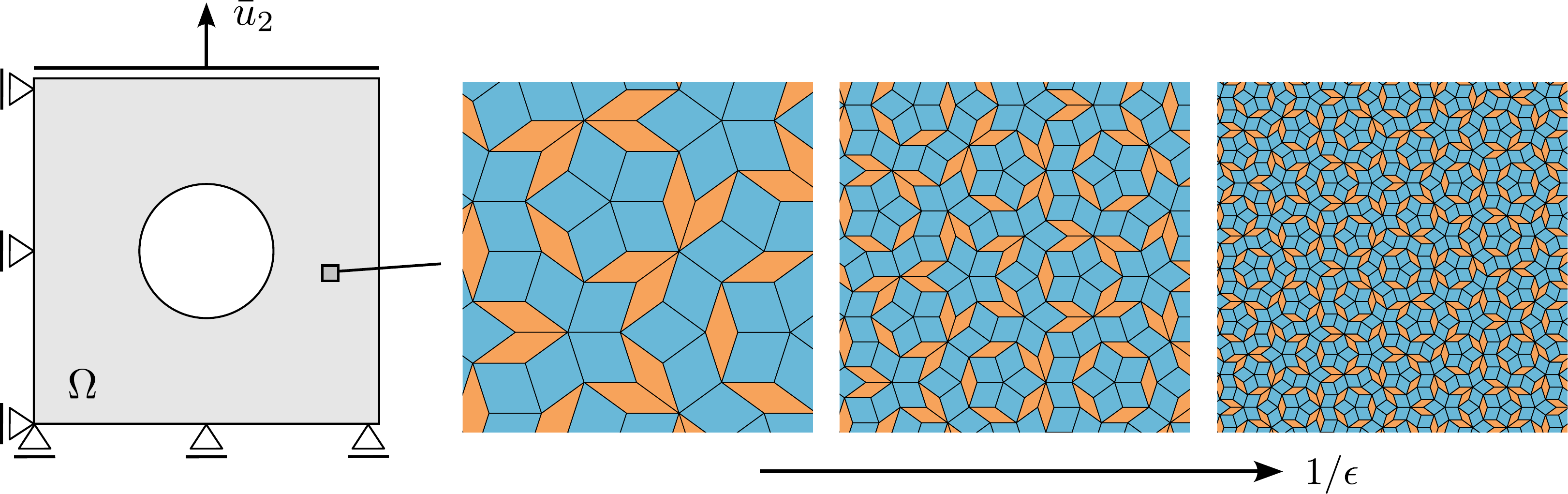}
    \caption{Boundary-value problem for the perforated plate with a Penrose tiling microstructure. The square plate has a length of 10~mm, and the central hole has a diameter of 4~mm. The specimen is under plane-stress conditions. For the direct numerical simulations, the microstructural length $\epsilon$ of the thick and thin rhombi (blue and orange, respectively) is taken progressively smaller while the macroscopic domain remains fixed.}
    \label{fig-penrose-domain}
\end{figure}

With the homogenized elasticity tensor determined, we investigate the elastic response of the two-dimensional Penrose material in a macroscopic boundary-value problem. The physical domain is a square of length $L=10$~mm with a circular hole of radius $r=2$~mm at its center (Figure~\ref{fig-penrose-domain}). The thick and thin rhombi are assigned Young's moduli
$E^{(1)}=8$~MPa and $E^{(2)}=1$~MPa, 
with a common Poisson's ratio $\nu=0.30$. We assume plane-stress conditions. Letting $x_1$ and $x_2$ denote the horizontal and vertical coordinates, we prescribe $\bar{u}_2=0.1$~mm on $x_2=L$, while setting $\bar u_2=0$ on $x_2=0$ and $\bar u_1=0$ on $x_1=0$.

\begin{figure}[t!]
    \centering
    \includegraphics[width=\linewidth]{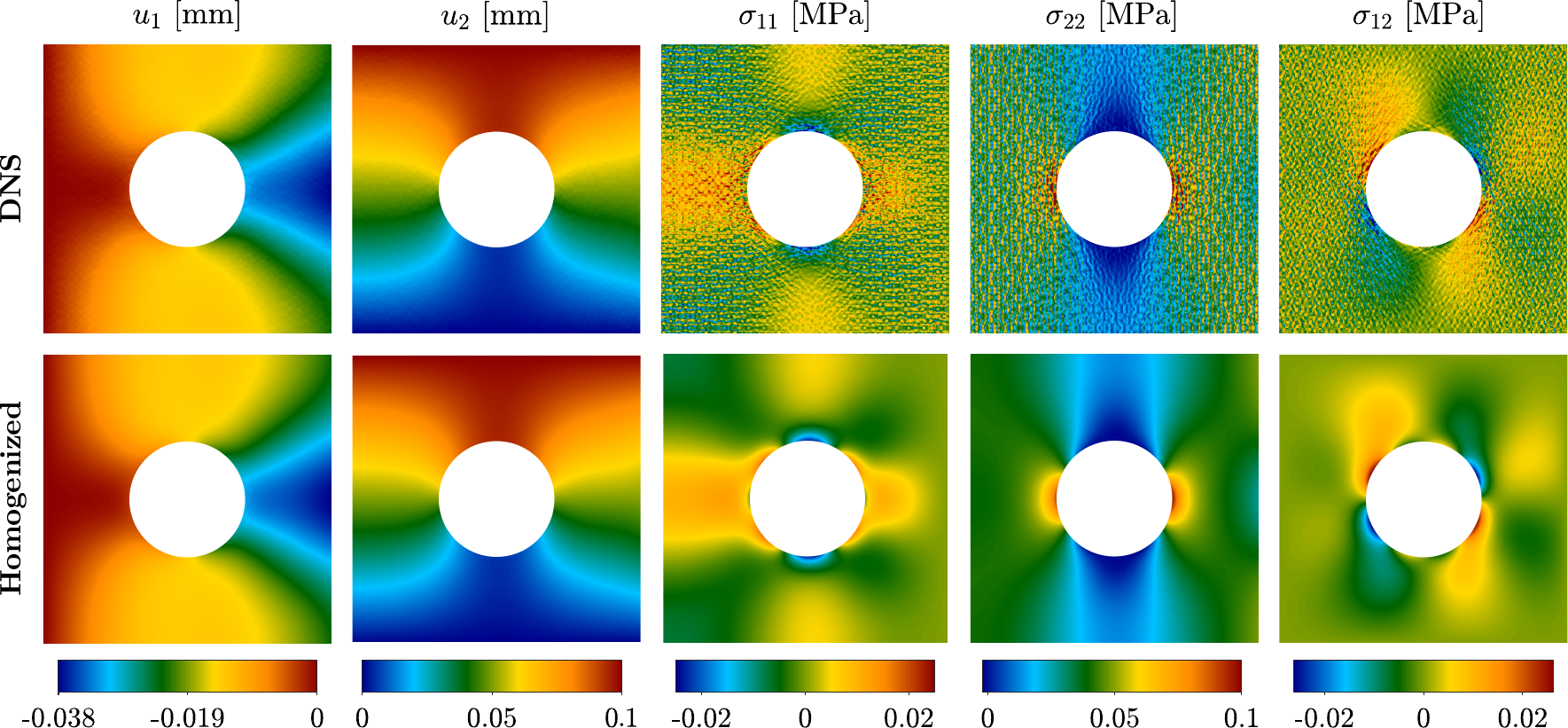}
    \caption{Comparison of displacement fields and stress fields obtained from direct finite element simulations (top) and the homogenized model (bottom) for the Penrose tiling at $\epsilon=0.01$.}
    \label{fig:Penrose_fields_str}
\end{figure}

We compare the response of the homogenized continuum model with direct numerical simulations (DNS) that explicitly resolve the material distribution in the Penrose tiling. Quadratic triangular elements are used for both direct and homogenized models. Mesh convergence studies were performed for both. Figure~\ref{fig:Penrose_fields_str} compares the displacement and stress fields obtained from one realization of the DNS at $\epsilon=0.01$ and the homogenized model. Agreement is observed for all quantities over the entire domain, indicating that the homogenized model captures the macroscopic response of the Penrose material. As expected, the stress fields obtained from the direct simulations exhibit local fluctuations associated with the heterogeneous microstructure, whereas the homogenized model produces smooth fields governed by the effective elastic properties. Nevertheless, both models reflect the influence of the circular hole on the deformation and stress distributions, with quantitative agreement.

\begin{figure}[b!]
    \vspace*{1em}
    \begin{center}
    \begin{tabular}{cc}
    \begin{overpic}[width=0.46\textwidth]{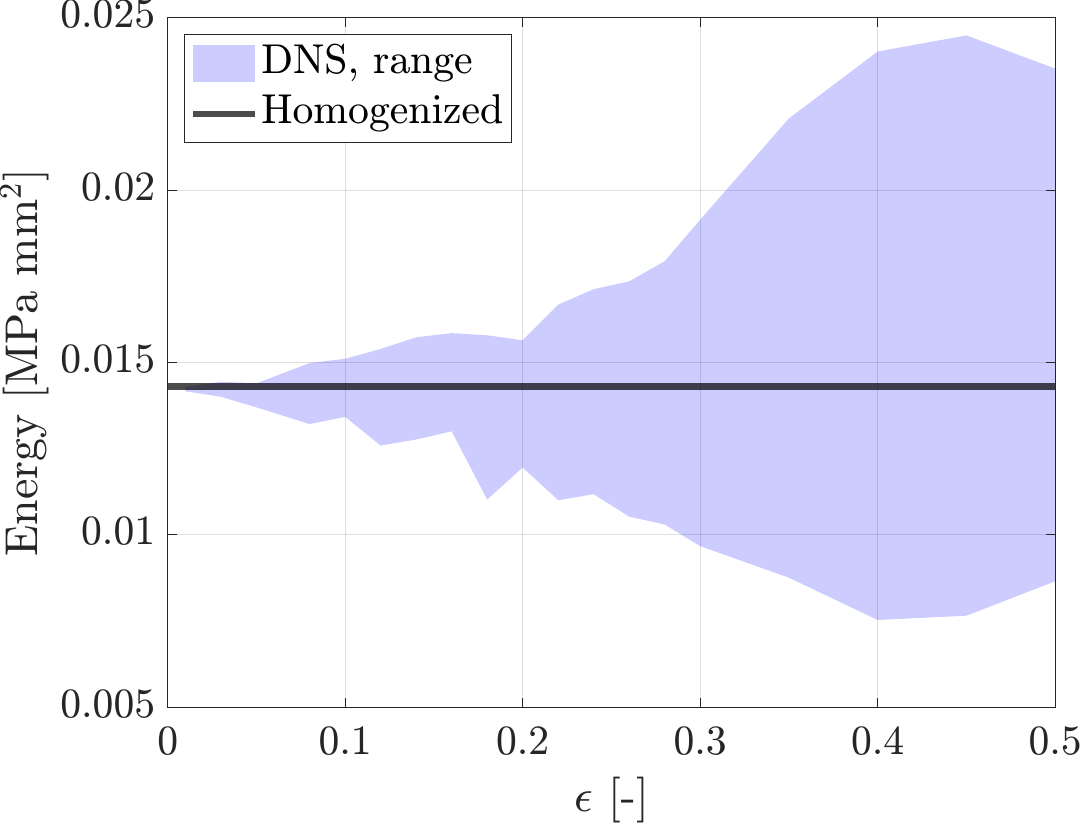}
        \put(0,2){\small{(a)}}
    \end{overpic} &
    \begin{overpic}[width=0.46\textwidth]{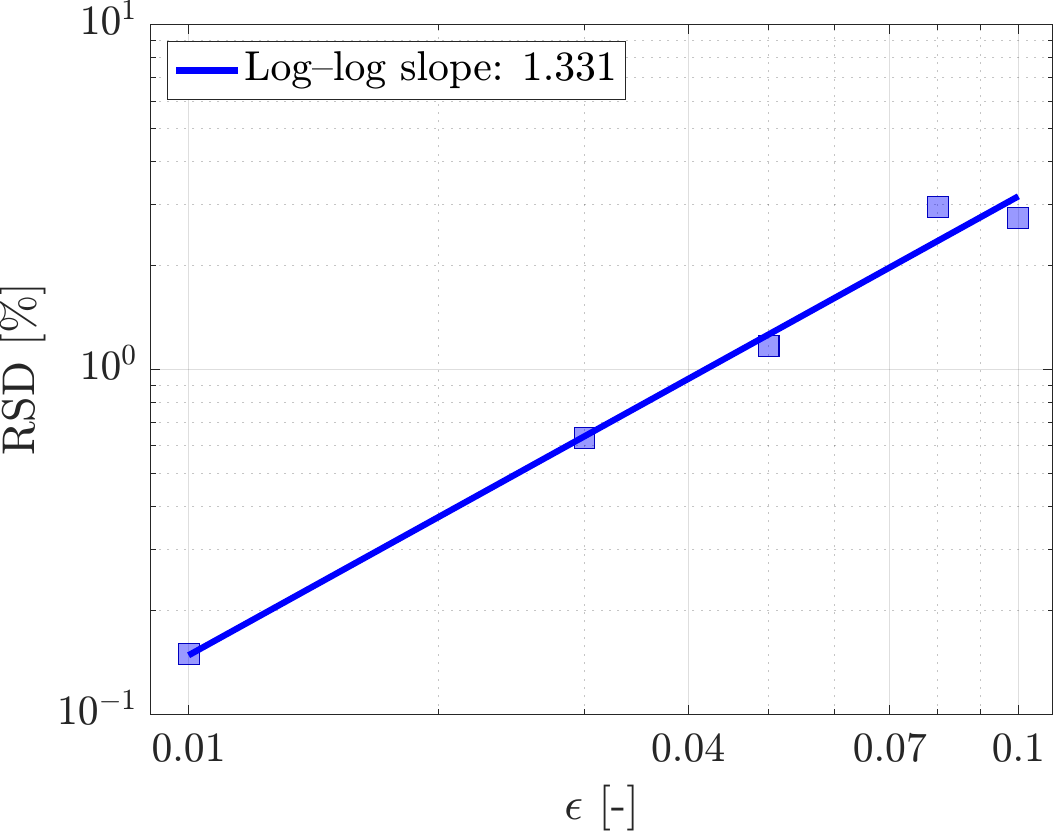}
        \put(0,2){\small{(b)}}
    \end{overpic}
    \end{tabular}
    \end{center}
    \vspace*{-1em}
    \caption{Energy convergence as $\epsilon\to 0$ for the Penrose tiling. (a) Energy range from direct numerical simulations (DNS) of 100 finite realizations (shaded region) and the homogenized limit (black line). (b) Decay of energy relative standard deviation (RSD) across realizations of the direct numerical simulations, showing a well-defined convergence rate.}
\label{fig:penrose_energy}
\end{figure}

We further examine the energy convergence, assessing the variability of the finite DNS patches as the characteristic microscopic length $\epsilon$ decreases. For each microstructural length scale $\epsilon$, $100$ phase-offset realizations are considered to assess the influence of the relative placement of the quasiperiodic structure. The first realization uses a prescribed offset vector $\bm\gamma_0$. For each of the remaining 99 realizations, indexed by $s=2,\ldots,100$, four independent random variables $\zeta_j^{(s)}$, $j=1,\ldots,4$, are drawn from the uniform distribution $U(0,1)$. These variables are added to the first four components according to $\gamma_j^{(s)}=\gamma_{0j}+\zeta_j^{(s)}$. The fifth component is chosen to enforce $\sum_{j=1}^{5}\gamma_j^{(s)}=0$.   

Figure~\ref{fig:penrose_energy}(a) shows the convergence of the strain energy. As $\epsilon$ decreases from $0.50$ to $0.01$, the range of energies over 100 DNS realizations progressively narrows and approaches the homogenized value. This behavior indicates that the influence of the phase offset diminishes with decreasing microstructural scale and provides numerical evidence of convergence towards the homogenized elastic response. Figure~\ref{fig:penrose_energy}(b) further shows that the relative standard deviation (RSD) of the strain energy across realizations decreases with $\epsilon$, with an estimated convergence rate of $1.331$.

Of course, the homogenized model offers substantially greater computational efficiency than the direct finite element model. Based on the mesh convergence study, the homogenized model requires only $9{,}908$ degrees of freedom to achieve a relative energy error
of $0.12\%$ compared with the extrapolated value, whereas the DNS at $\epsilon = 0.01$ requires $176{,}772$ degrees of freedom.

\subsection{Three-dimensional icosahedral tiling}

\begin{figure}[b!]
    \vspace*{1em}
    \begin{center}
    \begin{tabular}{cc}
    \begin{overpic}[width=0.43\textwidth]{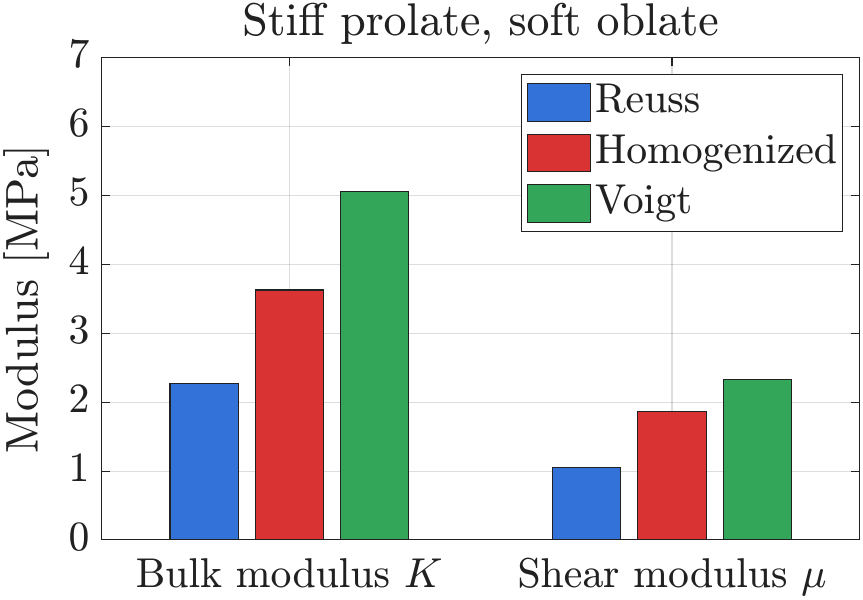}
        \put(0,2){\small{(a)}}
    \end{overpic} &
    \begin{overpic}[width=0.43\textwidth]{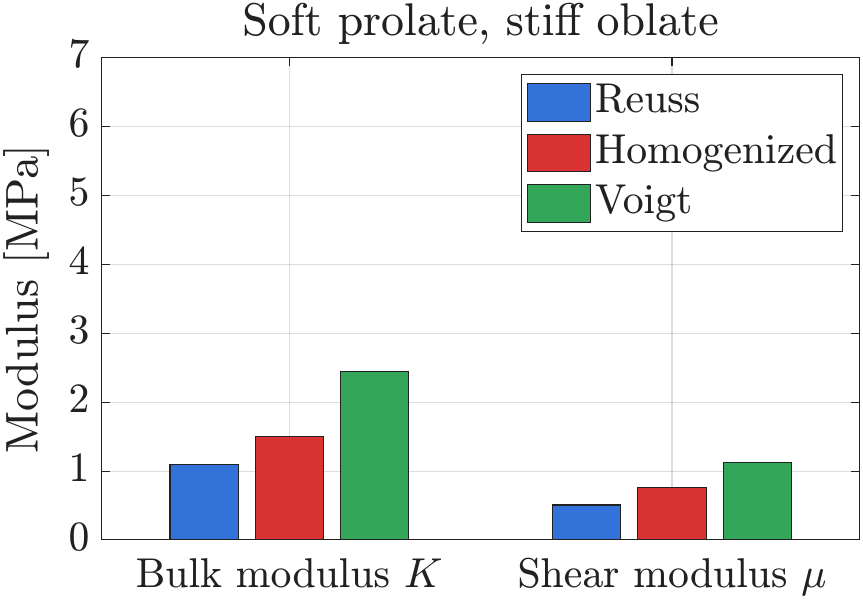}
        \put(0,2){\small{(b)}}
    \end{overpic}
    \end{tabular}
    \end{center}
    \vspace*{-1em}
    \caption{Comparison of the homogenized bulk and shear moduli
    with the Voigt and Reuss bounds for the AKN tiling with different material assignments.
    (a) Stiff prolate rhombohedra and soft oblate rhombohedra.
    (b) Soft prolate rhombohedra and stiff oblate rhombohedra.}
    \label{fig:akn_voigt_reuss_bounds}
\end{figure}

We compute the homogenized elasticity tensor for the AKN material by solving the projected cell problems with the FFT-based method using $N=31$ grid points along each of the six Fourier directions. We assign Young's moduli $E^{(1)}=8$~MPa and $E^{(2)}=1$~MPa to two types of rhombohedra, with a common Poisson's ratio $\nu=0.30$. By way of illustration, when we assign the stiff material to the prolate rhombohedra and the soft
material to the oblate rhombohedra, the resulting homogenized elasticity tensor is
\begin{equation}\label{eq:akn_hom_tensor}
    [\bm{\mathsf C}^{\rm eng}_0]=\begin{pmatrix}
    6.107&2.388&2.388&0&0&0\\
    2.388&6.107&2.388&0&0&0\\
    2.388&2.388&6.106&0&0&0\\
    0&0&0&1.859&0&0\\
    0&0&0&0&1.859&0\\
    0&0&0&0&0&1.859
    \end{pmatrix}\,\mathrm{MPa}.
\end{equation}
The tensor exhibits an approximately isotropic response, consistent with the icosahedral symmetry of the AKN tiling. In particular, the normal stiffness components are nearly equal, and the isotropic relation between normal and shear stiffness is closely satisfied up to numerical tolerance. 

Figure~\ref{fig:akn_voigt_reuss_bounds} compares the homogenized bulk and shear moduli of the AKN tiling with the Voigt and Reuss bounds for two different material assignments. The bounds are evaluated using the volume fractions of the prolate ($f_{\mathrm{prolate}}=\tau^2/(1+\tau^2)\approx72.36\%$) and oblate ($f_{\mathrm{oblate}}=1/(1+\tau^2)\approx27.64\%$) rhombohedra in the infinite AKN tiling. When the prolate rhombohedra are assigned the stiffer material, the homogenized bulk and shear moduli follow as $K_0=3.628$~MPa and $\mu_0=1.859$~MPa, respectively. These values lie between the corresponding Reuss bounds, $K_{\mathrm{R}}=2.272$~MPa and $\mu_{\mathrm{R}}=1.048$~MPa, and Voigt bounds, $K_{\mathrm{V}}=5.054$~MPa and $\mu_{\mathrm{V}}=2.333$~MPa. Inverting the phase assignment yields $K_0=1.505$~MPa and $\mu_0=0.757$~MPa, which similarly lie between the corresponding Reuss bounds, $K_{\mathrm{R}}=1.099$~MPa and $\mu_{\mathrm{R}}=0.507$~MPa, and Voigt bounds, $K_{\mathrm{V}}=2.446$~MPa and $\mu_{\mathrm{V}}=1.129$~MPa.

\begin{figure}[t!]
\centering
\includegraphics[width=\linewidth]{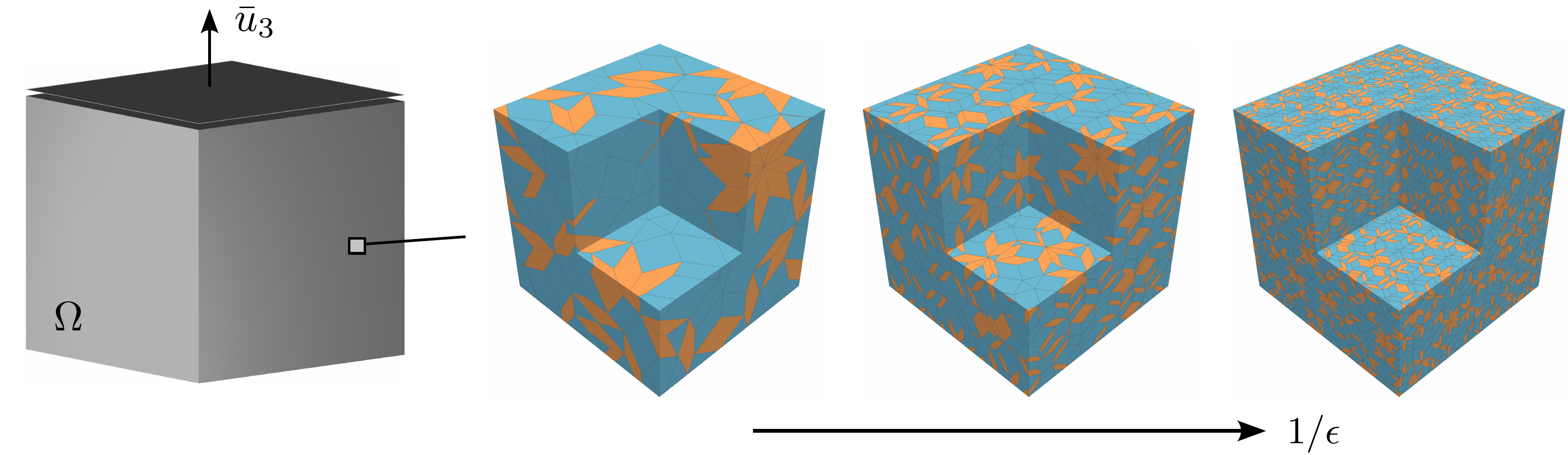}
\caption{Boundary-value problem for the cube with an AKN tiling microstructure. The cube has a length of 10~mm. Normal displacements are fixed on three adjacent faces, leaving the others free. For the direct numerical simulations, the microstructural length $\epsilon$ of the prolate and oblate rhombohedra (blue and orange, respectively) is taken progressively smaller while the macroscopic domain remains fixed (a corner is cut away to reveal the internal arrangement of the microstructure).}
\label{fig:icosahedral_cell}
\end{figure}

With the homogenized elasticity tensor determined, we investigate the elastic response of the three-dimensional AKN material in a macroscopic boundary-value problem. The physical domain is a cube of length $L=10$~mm (Figure~\ref{fig:icosahedral_cell}). The prolate and oblate rhombohedra are assigned Young's moduli $E^{(1)}=8$~MPa and $E^{(2)}=1$~MPa, respectively, with a common Poisson's ratio $\nu=0.30$. Let $x_1$, $x_2$, and $x_3$ denote the Cartesian coordinates, with $x_3$ oriented vertically. We set $\bar u_1=0$ on $x_1=0$, $\bar u_2=0$ on $x_2=0$, and $\bar u_3=0$ on $x_3=0$, while prescribing $\bar u_3=0.1$~mm on the top face $x_3=L$.

\begin{figure}[b!]
    \centering
    \includegraphics[width=0.85\linewidth]{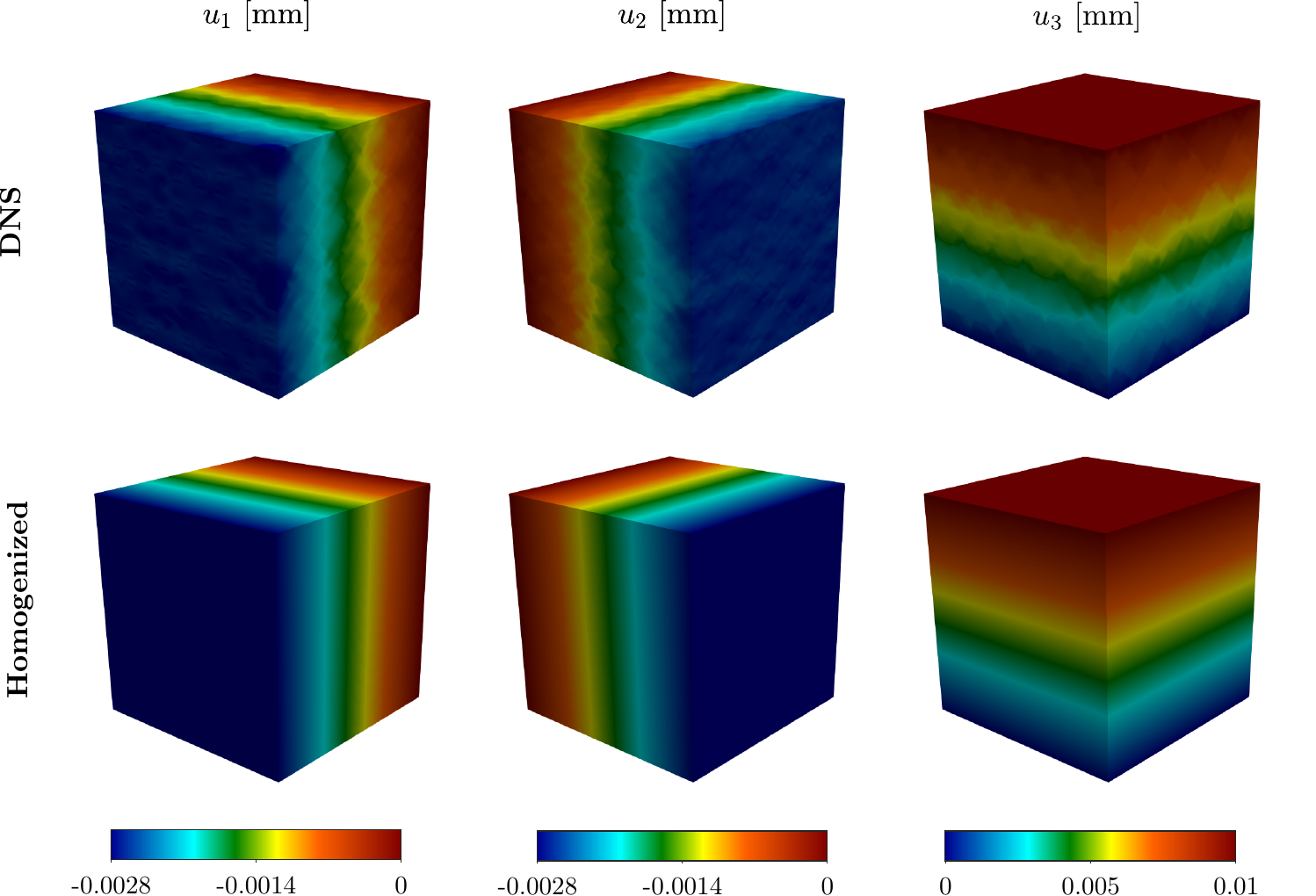}
    \caption{Comparison of displacement fields obtained from direct finite element simulations (top) and the homogenized model (bottom) for the AKN tiling at $\epsilon=0.05$.}
    \label{fig:ico_fields_kin}
\end{figure}

We compare the response of the homogenized continuum model with direct numerical simulations that explicitly resolve the material distribution in the AKN tiling. Linear tetrahedral elements are used for both the direct and homogenized models. Figure~\ref{fig:ico_fields_kin} compares the displacement fields obtained from one DNS realization at $\epsilon=0.05$ and the homogenized model. Despite the local fluctuations induced by the quasiperiodic microstructure, the DNS fields closely agree with those predicted by the homogenized model throughout the domain. In particular, the transverse displacement fields $u_1$ and $u_2$ exhibit similar spatial distributions and magnitudes, consistent with the effective Poisson's ratio of the quasiperiodic material.

\begin{figure}[t!]
    \vspace*{1em}
    \begin{center}
    \begin{tabular}{cc}
    \begin{overpic}[width=0.46\textwidth]{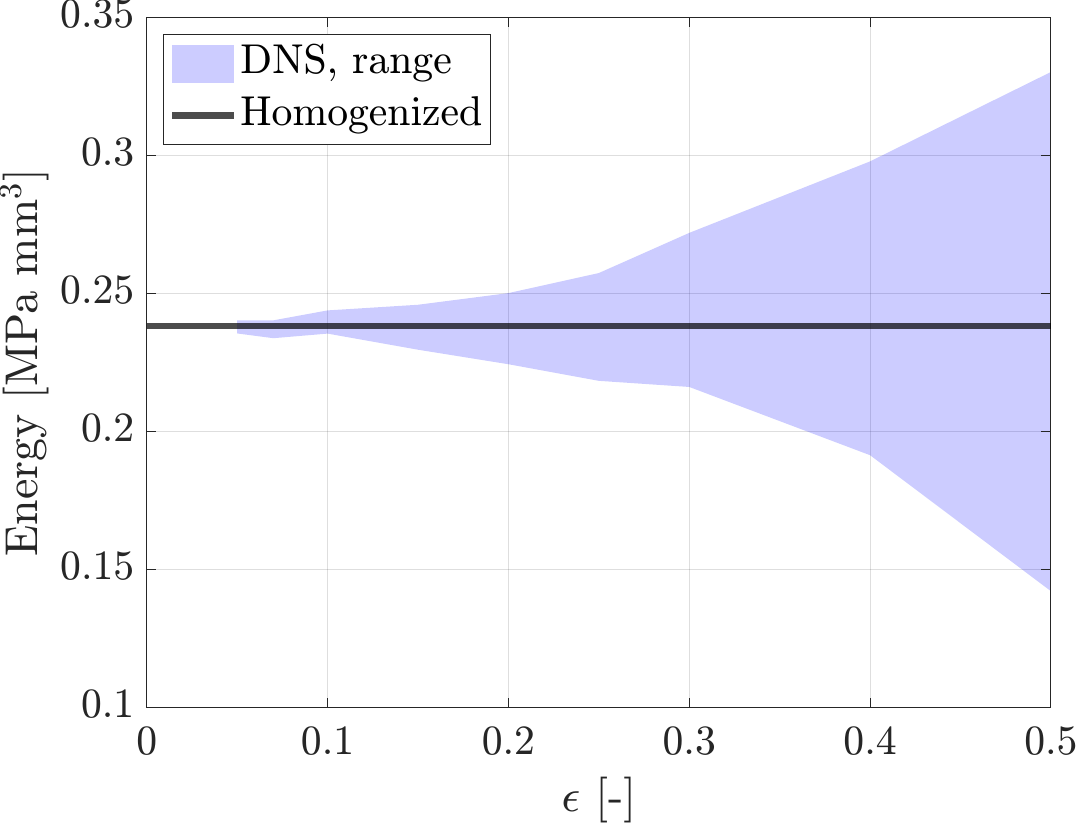}
        \put(0,2){\small{(a)}}
    \end{overpic} &
    \begin{overpic}[width=0.46\textwidth]{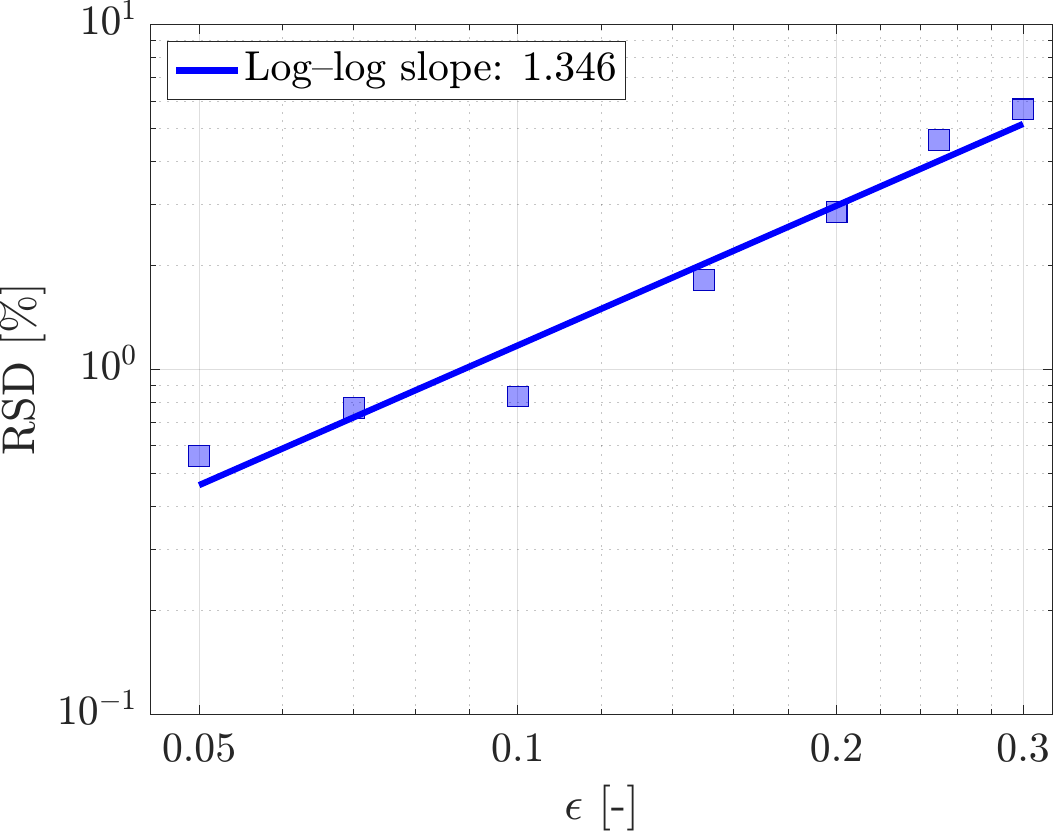}
        \put(0,2){\small{(b)}}
    \end{overpic}
    \end{tabular}
    \end{center}
    \vspace*{-1em}
    \caption{Energy convergence as $\epsilon\to 0$ for the icosahedral tiling. (a) Energy range from DNS of 30 finite realizations (shaded region) and the homogenized limit (black line). (b) Decay of energy RSD across realizations of the direct numerical simulations, showing a well-defined convergence rate.}
\label{fig:icosahedral_energy}
\end{figure}

We further examine the energy convergence and the variability across finite DNS realizations as the characteristic microstructural length $\epsilon$ decreases. For each value of $\epsilon$, $30$ phase-offset realizations are considered to assess the influence of the relative placement of the quasiperiodic structure. Figure~\ref{fig:icosahedral_energy}(a) shows that, as $\epsilon$ decreases from $0.50$ to $0.05$, the range of energies across these realizations progressively narrows and approaches the homogenized value. This behavior provides numerical evidence that the three-dimensional elastic response of the AKN material converges toward its homogenized limit. Figure~\ref{fig:icosahedral_energy}(b) further shows that the relative standard deviation of the strain energy decreases with $\epsilon$, with an estimated convergence rate of $1.346$.

\section{Conclusions and perspectives}
\label{sec:conclusions}

We have developed a computational framework for the homogenization of quasiperiodic elastic materials based on their higher-dimensional periodic representations. By exploiting the cut-and-projection structure of the quasiperiodic materials, the absence of a periodic unit cell in physical space is overcome through the formulation of periodic cell problems in higher-dimensional spaces. These problems are solved using an FFT-based approach to determine the effective elasticity tensors, including symmetry-preserving discretization schemes. The framework has been implemented for two-dimensional Penrose tilings and three-dimensional AKN icosahedral tilings.

The computed effective elasticity tensors exhibit approximately isotropic responses, consistent with the rotational symmetries of the quasiperiodic architectures. To assess the homogenized models, direct finite element simulations resolving the quasiperiodic microstructures were performed over progressively decreasing characteristic tile sizes. In both two and three dimensions, the direct simulations approach the homogenized energies as the microstructural scale decreases. Moreover, the variability across different phase realizations diminishes, indicating that the influence of the relative placement of the microstructure becomes negligible in the homogenization limit. These results demonstrate the applicability of the homogenization of quasiperiodic elastic materials and provide a computational approach for evaluating their effective properties. The accompanying open-source implementations provide a basis for investigating other quasiperiodic architectures and their effective mechanical behavior.

Several directions for future research arise from the present work. A natural extension is the consideration of nonlinear mechanical behavior, including elastoplasticity and buckling phenomena, for which the interaction between quasiperiodic geometry and local deformation mechanisms can play an important role. Furthermore, the proposed framework provides a foundation for addressing inverse homogenization problems in quasiperiodic materials, where the underlying microstructure is designed to achieve prescribed effective properties. Additionally, the proposed framework could be extended to discrete-to-continuum homogenization, enabling the derivation of effective continuum models for quasicrystalline lattices.

\bibliography{biblio}
\end{document}